\input amstex
\mag = \magstep1
\documentstyle{amsppt}
\document
%\newsymbol \rtimes 226F
\TagsOnRight
\NoBlackBoxes

\define\a{{r}}
\define\abel{^{\text{abel}}}
\define\An{{\operatorname{Aut}(F_n)}}
\define\Anbar{{\overline{A_n}}}
\define\Aut{{\operatorname{Aut}}}
\define\bb{{\varsigma}}
\define\bC{\Bbb C}
\define\bZ{\Bbb Z}
\define\bQ{\Bbb Q}

\define\bT{\Bbb T}
\define\Coker{{\operatorname{Coker}}}
\define\diag{{\operatorname{diag}}}
\define\ebar{{\overline{e}}}
\define\Fn{{F_n}}
\define\GL{{\operatorname{GL}}}
\define\h{{\overline{h}}}
\define\Hom{{\operatorname{Hom}}}
\define\HZ{{H_{\Bbb Z}}}
\define\IA{{\operatorname{IA}}}
\define\IAn{{IA_n}}
\define\Inn{{\operatorname{Inn}(F_n)}}
\define\inv{^{-1}}
\define\IOn{{IO_n}}
\define\Ker{{\operatorname{Ker}}}
\define\La{{\Lambda}}
\define\M1g{{\Cal M_{g, 1}}}
\define\Mbar{{\overline{\Cal M}}}
\define\Mzbar{{\overline{\Cal M^0}}}
\define\Mg{{\Cal M_{g}}}
\define\Mgstar{{\Cal M_{g, *}}}
\define\On{{\operatorname{Out}(F_n)}}
\define\R{{\operatorname{R}}}
\define\RFn{{\widehat{R[F_n]}}}
\define\tensors{{{H_{\Bbb Z}}^*\otimes\Lambda^2H_{\Bbb Z}}}
\define\SL{{\operatorname{SL}}}
\define\sgn{\operatorname{sgn}}
\define\std{{\operatorname{std}}}
\define\T{{\widehat{T}}}

\topmatter
\title 
Cohomological Aspects of Magnus Expansions.
\endtitle
\dedicatory
Dedicated to Professor Tatsuo Suwa on his sixtieth birthday
\enddedicatory
\date{September 7, 2006}\enddate
\rightheadtext{Cohomological Aspects of Magnus Expansions}
\author Nariya KAWAZUMI
\endauthor
\affil
Department of Mathematical Sciences,\\ University of Tokyo
\endaffil
\address
Tokyo, 153-8914 Japan 
\endaddress
\email
kawazumi\@ms.u-tokyo.ac.jp
\endemail

\abstract
We generalize the notion of a Magnus expansion of a free group
in order to extend each of the Johnson homomorphisms 
defined on a decreasing filtration of the Torelli
group for a surface with one boundary component
to the whole of the automorphism group of a free group 
 $\operatorname{Aut}(F_{n})$.
The extended ones are {\it not} homomorphisms, but satisfy
an infinite sequence of coboundary relations, so that
we call them {\it the Johnson maps}.
In this paper we confine ourselves to studying the first and 
the second relations, which have cohomological consequences 
about the group $\operatorname{Aut}(F_{n})$ and 
the mapping class groups for surfaces. 
The first one means that the first Johnson map is
a twisted $1$-cocycle of the group $\operatorname{Aut}(F_{n})$.
Its cohomology class coincides with ``the unique elementary 
particle" of all the Morita-Mumford classes on the mapping 
class group for a surface [Ka1] [KM1]. 
The second one restricted to the mapping class group is equal to a
fundamental relation among twisted Morita-Mumford classes proposed by
Garoufalidis and Nakamura [GN] and established
 by Morita and the author [KM2]. 
This means we give a simple and coherent proof 
of the fundamental relation. 
The first Johnson map gives the abelianization 
of the induced automorphism group
$IA_n$ of a free group in an explicit way.
\endabstract

\subjclass{Primary 20J05.
  Secondary 14H10, 20E05, 20F28, 57R20}\endsubjclass
\keywords{Magnus expansion, Johnson homomorphisms,
Morita-Mumford class, automorphism group of a free group, 
mapping class group, IA-automorphism group}\endkeywords
\endtopmatter

\head Introduction\endhead

In the cohomological study of the mapping class group for a surface, 
or equivalently, in that of the moduli space of Riemann surfaces, 
the Morita-Mumford class, $e_m = (-1)^{m+1}\kappa_m$, $m \geq 1$, 
[Mu] [Mo1] plays an essential role. 
Harer [Har] established that the $i$-th cohomology group of 
the mapping class group of genus $g$ does not depend 
on the genus $g$ if $i < g/3$. 
This enables us to consider the stable cohomology algebra 
of the mapping class groups, $H^*(\Cal M_\infty; \Bbb Z)$. 
Recently Madsen and Weiss [MW] proved that the rational stable 
cohomology algebra of the mapping class groups,
$H^*(\Cal M_\infty; \Bbb Q)$, 
is generated by the Morita-Mumford classes. \par

As was shown by Miller [Mi], each of the Morita-Mumford classes is 
indecomposable in the algebra $H^*(\Cal M_\infty; \Bbb Q)$. 
But, if we consider the twisted cohomology of the mapping 
class groups, then the Morita-Mumford class decomposes itself 
into an algebraic combination of some copies of the extended 
first Johnson homomorphism $\tilde k$, introduced by Morita [Mo3], 
by means of the intersection product of the surface [Mo4] [KM1] [KM2].  
In short, ``the unique elementary particle" for all the
Morita-Mumford classes is the extended first Johnson homomorphism 
$\tilde k$. \par

Throughout this  paper, we often confine ourselves to the
mapping class group for a compact surface of genus $g$ with $1$ boundary
component, $\Cal M_{g, 1}$, for simplicity.
From Harer's stability theorem [Har] the $i$-th integral cohomology 
group of the group $\Cal M_{g, 1}$ is isomorphic to 
$H^*(\Cal M_\infty; \Bbb Z)$ if $i < g/3$. 
The extended first Johnson homomorphism $\tilde k$ 
is a crossed homomorphism defined on
$\Cal M_{g, 1}$ [Mo3] and equals to $\frac16 m_{0, 3}$,  
which is the $(0, 3)$-twisted (or generalized) Morita-Mumford 
class [KM2] in a terminology of the author's work [Ka1]. 
There are many nontrivial relations among different
combinations of $\tilde k$, 
so that all the cohomology classes
with trivial coefficients we obtain from $\tilde k$ are just 
the polynomials of the Morita-Mumford classes [KM1] [KM2]. 
Here we should remark our sign convention on cap products 
in the present paper is different from our previous one 
in [KM1] [KM2], where $\tilde k$ equals to $-\frac16 m_{0, 3}$. 
It will be explained in \S5. \par

Let $H$ denote the first integral homology group for the surface, and 
$H^*$ its dual, namely, the first integral cohomology group of 
the surface. The twisted Morita-Mumford class $m_{i, j}$ is a twisted 
cohomology class in $H^{2i+j-2}(\Cal M_{g, 1}; \La^jH)$ [Ka1]. 
When $j=0$, it equals to the original Morita-Mumford class $e_i = 
m_{i+1, 0}$. In his study of Riemann constants 
on Jacobi varieties Earle [E], p\.272, had already introduced 
an integral $1$-cocycle $(2-2g)\psi$ of $\Cal M_{g, 1}$ 
with coefficients in $H$.  It represents the $(1, 1)$-twisted
Morita-Mumford class [Mo1], p\.81, [Ka1], p\.147. 
From a similar reason to the trivial coefficients, 
all the cohomology classes with twisted coefficients 
we obtain from $\tilde k$ are just the polynomials of 
the twisted Morita-Mumford classes [KM1] [KM2]. \par
  
The natural action of the mapping class group 
on the fundamental group of the compact surface, 
which is isomorphic to a free group of rank $2g$, 
defines a homomorphism $\M1g \to \operatorname{Aut}(F_{2g})$.
Here $F_n$ is a free group of rank $n \geq 2$.
Andreadakis [An] introduced a decreasing filtration of the group $\An$ 
by using the action of $\An$ on the lower central series of the group $\Fn$. 
Its pullback to the mapping class group $\M1g$ coincides with 
the decreasing filtration $\Cal M_{g, 1}(p)$, $p \geq 0$, 
introduced by Johnson [J2]. Moreover Johnson [J1] [J2] defined 
a sequence of injective homomorphisms, 
$\tau_p: \Cal M_{g, 1}(p)/\Cal M_{g, 1}(p+1) \to 
H^*\otimes \Cal L_{p+1}$, $p\geq 1$. 
The homomorphism $\tau_p$ is called the $p$-th Johnson homomorphism. 
Here $\Cal L_p$ denotes the $p$-th component 
of the free Lie algebra generated by $H$.  
The group $\Cal M_{g, 1}(0)$ is the full mapping class group 
$\Cal M_{g, 1}$, and $\Cal M_{g, 1}(1)$ the Torelli group 
$\Cal I_{g, 1}$.
Johnson [J3] proved $\tau_1$ induces a surjection 
${\Cal I_{g,1}}^{\text{abel}}\to \Lambda^3H$, 
and that its kernel is $2$-torsion. 
The extended one $\tilde k$ coincides with
$\tau_1$ on the Torelli group $\Cal M_{g, 1}(1) = \Cal I_{g, 1}$.\par 

The purpose of the present paper is to give a coherent point of view 
about all the Johnson homomorphisms, the twisted Morita-Mumford classes 
and their higher relations. 
It should be remarked that Hain [Hai] established 
a remarkable theory about all the Johnson homomorphisms 
on the Torelli groups by means of Hodge theory. 
On the other hand we focus our study on the automorphism group
$\operatorname{Aut}(F_{n})$ 
instead of the mapping class group $\Cal M_{g, 1}$, and so 
all the tools we will use are contained in elementary algebra, 
except the Lyndon-Hochschild-Serre spectral sequences 
in group cohomology [HS].
As a consequence all the Johnson homomorphisms extend 
to the whole automorphism group $\operatorname{Aut}(F_{n})$ in a natural way.
But they are {\it not} homomorphisms any longer. 
They satisfy an infinite sequence of coboundary relations, 
so that we call them {\it the Johnson maps}.
In a forthcoming paper [Ka3] we will show 
that how far they are from true homomorphisms is measured by the Stasheff
associahedrons [S] in an infinitesimal way. \par 

The key to our consideration is the notion of the Magnus expansion 
of a free group. 
In the classical context it is an embedding of the free group into the
multiplicative group of the completed tensor algebra generated by the
abelianization of the free group, ${F_n}^{\text{abel}}$.
Magnus [M1] constructed it by an explicit use of the generators. 
It can be defined also by using Fox' free differentials. 
See, e.g., [F].   
Throughout this paper we call the classical one {\it the standard
Magnus expansion}. 
It is remarkable that Kitano [Ki] described explicitly 
the Johnson homomorphism $\tau_p$ on the group $\Cal M_{g, 1}(p)$ 
in terms of the standard Magnus expansion. 
On the other hand, Bourbaki [Bou] gave a general theory 
on an embedding of the free group into the multiplicative group 
of the completed tensor algebra. 
Our consideration is a refinement of Kitano's
description and Bourbaki's theory. 
More precisely, we generalize the notion of a Magnus expansion. 
We call a map of the free group $F_n$ into the completed
tensor algebra  with a property enough to define the Johnson homomorphisms {\it
a Magnus expansion} (\S1). 
The set of all the Magnus expansions in our
generalized sense, $\Theta_n$, admits two kinds of group actions.  
The automorphism group $\operatorname{Aut}(F_{n})$ and 
a certain Lie group $IA(\widehat{T})$ act on it in a natural way. 
The latter action is free and transitive. 
So, if we fix a Magnus expansion $\theta\in \Theta_n$, 
then we can write the action of $\operatorname{Aut}(F_{n})$ on the set
$\Theta_n$ in terms of the action of $IA(\widehat{T})$. 
This gives us the $p$-th Johnson map
$\tau^\theta_p: \operatorname{Aut}(F_{n}) \to H^*\otimes H^{\otimes p+1}$ 
for each $p \geq 1$, which is equal to the original  $\tau_p$ on the subgroup
$\Cal M_{g, 1}(p)$. 
On the mapping class group $\Cal M_{g, 1}$ the first one
$\tau^\theta_1$ coincides with $\tilde k$, and $\tau^\theta_2$ gives a
fundamental relation among combinations of $\tilde k$'s.  
The cohomology class of $\tau^\theta_1$ is equal to the Gysin image of the
square of a certain twisted first cohomology class $k_0$ on the semi-direct
product $F_n\rtimes\operatorname{Aut}(F_{n})$.
The second Johnson map $\tau^\theta_2$ gives a fundamental relation among
twisted Morita-Mumford classes proposed by Garoufalidis and Nakamura [GN] 
and established by Morita and the author [KM2]. 
Recently Akazawa [Ak] gave an alternative proof of it 
by using representation theory of the symplectic group. 
This means we give a simple and coherent proof 
of the fundamental relation (Theorem 5.5). 
In \S 5 we show some of the twisted Morita-Mumford classes, 
$m_{0, j}$ and $m_{1, j}$,  
extend to the automorphism group $\operatorname{Aut}(F_{n})$. \par 

In a similar way to the Torelli group $\Cal I_{g, 1}$ 
the IA-automorphism group $IA_n$ is defined to be the kernel 
of the homomorphism $\operatorname{Aut}(F_{n}) \to GL_n(\Bbb Z)$ 
induced by the natural action on the abelianization $H = {F_n}^{\text{abel}}
\cong \Bbb Z^n$. 
In \S 6 we prove the map $\tau^\theta_1$ induces an isomorphism
${IA_n}^{\text{abel}} \cong H^*\otimes \Lambda^2 H \cong \Bbb Z^{\frac12
n^2(n-1)}$ using Magnus' generators of the group $IA_n$ [M2].\par 

Some of the results in this paper were
announced in [Ka2]. In a forthcoming paper [Ka3] we will study  
the Teichm\"uller space $\Cal T_{g, 1}$ of triples $(C, P_0, v)$, 
where $C$ is a compact Riemann surface of genus $g$, $P_0 \in C$, 
and $v \in T_{P_0}C\setminus\{0\}$. 
There we will define and study a canonical map of
$\Cal T_{g, 1}$ into $\Theta_{2g}$, 
which will make the advantage of our
generalization of Magnus expansions clearer.\par
\smallskip\noindent
{\it Acknowledgements}: 
The author would like to thank Shigeyuki Morita, 
Hiroaki Nakamura, Atsushi Matsuo, Teruaki Kitano, 
Yuuki Tadokoro and Takao Satoh for helpful discussions, 
and Soren Galatius for giving the author some informations 
about his unpublished work [Ga]. 
\par

\demo{Contents}\par
Introduction.\par
\S1. Magnus Expansions.\par
\S2. Johnson Maps.\par
\S3. Lower Central Series.\par
\S4. Twisted Cohomology Classes.\par
\S5. Mapping Class Groups.\par
\S6. The abelianization of $\IAn$.\par
\S7. Decomposition of Cohomology Groups.\par
\enddemo

\beginsection 1. Magnus Expansions.\par

In our generalized sense we define the notion 
of a Magnus expansion of a free group of rank $n$, $\Fn$. 
It is defined to be a group homomorphism of the free group 
into the multiplicative group of the completed tensor algebra 
$\T$ generated by the first homology group of the free group 
with a certain condition (Definition 1.1). 
An automorphism group $\IA(\T)$ of the algebra $\T$ acts on 
the set $\Theta_n$ consisting of all the Magnus expansions 
as well as the automorphism group  $\Aut(\Fn)$ of the free group. 
In the latter half of this section we introduce and study 
the group $\IA(\T)$. 
We prove the group $\IA(\T)$ acts on the set $\Theta_n$ 
in a free and transitive way (Theorem 1.3). 
Our construction of Johnson maps is based 
on this free and transitive action.
\par

Let $n \geq 2$ be an integer, $\Fn$ a free group of rank $n$
with free basis $x_1, x_2, \dots, x_n$
$$
\Fn = \langle x_1, x_2, \dots, x_n\rangle,
$$
and $R$ a commutative ring with a unit element $1$. 
We denote by $H = H_R$ the first homology group of 
the free group $\Fn$ with coefficients in $R$
$$
H = H_R := H_1(\Fn; R) 
= \Fn\abel\otimes_\bZ R \cong R^{\oplus n}.
$$
Here $G\abel$ is the abelianization of a group $G$, 
$G\abel = G/[G, G]$. We denote 
$$
[\gamma] := (\gamma \bmod [\Fn, \Fn])\otimes_\bZ 1 \in H
$$
for $\gamma \in \Fn$, and $X_i := [x_i] \in H$ 
for $i$, $1 \leq i\leq n$. 
The set $\{X_1, X_2, \dots, X_n\}$ is an $R$-free 
basis of $H$. \par
For the rest of the paper we write simply $\otimes$ and 
$\Hom$ for the tensor product 
and the homomorphisms over the ring $R$, respectively. 
The completed tensor algebra generated by $H$ 
$$
\T = \T(H) : = {\prod}^\infty_{m=0} H^{\otimes m}
$$
is equal to the ring of noncommutative formal power series 
$R\left<\left<X_1, X_2, \dots, X_n\right>\right>$, 
which Bourbaki [Bou] denotes by 
$\widehat{A}(\{x_1, x_2, \dots, x_n\})$.
The two-sided ideals
$$
\T_p := {\prod}_{m\geq p} H^{\otimes m}, \quad p \geq 1,
$$
give a decreasing filtration of the algebra $\T$.
For each $m$ we regard $H^{\otimes m}$ as a subspace of $\T$ 
in an obvious way. So we can write 
$$
z = \sum^\infty_{m=0} z_m 
= z_0 + z_1 + z_2 + \cdots + z_m + \cdots
$$
for $z = (z_m) \in \T$, $z_m \in H^{\otimes m}$.
It should be remarked that the subset $1+\T_1$ is 
a subgroup of the multiplicative group of the algebra $\T$, which 
Bourbaki [Bou] denotes by $\Gamma(\{x_1, x_2, \dots, x_n\})$ 
and calls {\it the Magnus group} 
over the set $\{x_1, x_2, \dots, x_n\}$. 
Now we can define a Magnus expansion of the free group $\Fn$ 
in our generalized sense. 
\proclaim{Definition 1.1} A map $\theta: \Fn \to 1 + \T_1$ is 
an $R$-valued Magnus expansion of the free group $\Fn$,
if 
\roster
\item $\theta: \Fn \to 1 + \T_1$ is a group homomorphism, and
\item $\theta(\gamma) \equiv 1 + [\gamma] \pmod{\T_2}$ for any
$\gamma 
\in \Fn$.
\endroster
\endproclaim
We write $\theta(\gamma) = \sum^\infty_{m=0}\theta_m(\gamma)$, 
$\theta_m(\gamma) \in H^{\otimes m}$. The $m$-th component 
$\theta_m: \Fn \to H^{\otimes m}$ is a map, but {\it not} 
a group homomorphism. The condition (2) is equivalent
to  the two conditions $\theta_0(\gamma) = 1$ and 
$\theta_1(\gamma) = [\gamma]$ for any $\gamma \in \Fn$.\par
From the universal mapping property of the free group $\Fn$,
for any $\xi_i \in \T_2$, $1 \leq i \leq n$, 
there exists a unique Magnus expansion $\theta$ satisfying 
$\theta(x_i) = 1 + X_i + \xi_i$ for each $i$. 
In other words, when we denote by $\Theta_n = \Theta_{n, R}$ 
the set of all the $R$-valued Magnus expansions, 
we have a bijection
$$
\Theta_n \cong (\T_2)^n, \quad
\theta \mapsto (\theta(x_i) - 1- X_i). \tag 1.1
$$
The standard Magnus expansion Magnus [M1] introduced
corresponds to $(0, 0, \dots, 0) \in (\T_2)^n$, 
which we denote
$$
\std: \Fn \to 1 + \T_1, \quad x_i \mapsto 1+X_i.
$$
\par
We denote by $\Aut(\T)$ the group of all the
filtration-preserving $R$-algebra automorphisms 
of the algebra $\T$. 
Here an $R$-algebra automorophism $U$ of $\T$ 
is defined to be {\it filtration-preserving} 
if $U(\T_p) = \T_p$ for each $p \geq 1$. 
It is easy to see whether an $R$-algebra endomorphism of 
$\T$ is a filtration-preserving automorphism or not.

\proclaim{Lemma 1.2} An $R$-algebra 
endomorphism $U$ of $\T$ is a filtration-preserving 
$R$-algebra automorphism, $U \in \Aut(\T)$, if and only if 
it satisfies the conditions
\roster
\item"(i)" $U(\T_p) \subset \T_p$ for each $p \geq 1$, and 
\item"(ii)" the induced 
endomorphism $\vert U\vert$ of $ \T_1/\T_2 = H$ is an isomorphism.
\endroster
\endproclaim

\demo{Proof}
We may suppose the $R$-algebra endomorphism $U$ satisfies 
the condition (i). We have 
$\Ker\, U =\Ker\left(U\vert_{\T_1}\right)$ and
$\Coker\, U = \Coker\left(U\vert_{\T_1}\right)$, 
since the algebra endomorphism $U$ preserves the direct sum 
decomposition $\T = R\cdot 1 \oplus \T_1$.
Hence the homomorphism of short exact sequences
$$
\CD
0 @>>> \T_2 @>>> \T_1 @>>> H @>>> 0\\
 @. @V{U}VV @V{U}VV @V{\vert U\vert}VV @.\\
0 @>>> \T_2 @>>> \T_1 @>>> H @>>> 0\\
\endCD
$$
induces a long exact sequence
$$
\multline
0 \to \Ker\left(U\vert_{\T_2}\right) \to \Ker\, U \to \Ker\vert U\vert
\to \Coker\left(U\vert_{\T_2}\right) \\
\to \Coker\, U \to \Coker\vert U\vert \to 0,
\quad \text{(exact)}.
\endmultline
\tag 1.2
$$
Here we prove 
$$
\Coker\left(U\vert_{\T_2}\right) = 0\quad
\text{if $\Coker\vert U\vert = 0$.}
\tag 1.3
$$
It suffices to show that, under the condition $\vert U\vert$ is 
surjective, for any given $w = (w_m) \in \T_2$, 
there exists a solution
$z = (z_m) \in \T_2$, $z_m \in H^{\otimes m}$, for the equation
$$
Uz = w. \tag 1.4
$$
The $m$-th component is
$$
\split
& \vert U\vert^{\otimes 2} z_2 = w_2,\quad \text{for $m = 2$,}\\
& \vert U\vert^{\otimes m}(z_m)
+ \text{terms in $z_2, z_3, \dots, z_{m-1}$}
= w_m,
\quad \text{for $m \geq 3$}.
\endsplit
$$
For each $m$, $\vert U\vert^{\otimes m}$ is surjective.
Hence we can find a tensor $z_m$ satisfying the equation
by induction on degree $m$. This proves (1.3).
\par
Now suppose $U \in \Aut(\T)$. Then 
we have $\Coker\,\vert U\vert = 0$ from (1.2), and 
$\Coker\left(U\vert_{\T_2}\right) = 0$ from (1.3). 
The sequence (1.2) implies
$\Ker \vert U\vert = 0$.
Hence $\vert U\vert$ is an isomorphism, that is, 
$U$ satisfies the condition (ii). 
\par
Conversely suppose $U$ is an $R$-algebra endomorphism 
satisfying the conditions (i) and (ii). 
Using the sequence (1.2) and (1.3), 
we obtain $\Coker\, U = 0$. 
Let $z = (z_m)$ be an element of 
$\Ker\left(U\vert_{\T_2}\right)$. 
The equation $Uz = 0$ is equivalent to 
$$
\split
& \vert U\vert^{\otimes 2}z_2 = 0 \in H^{\otimes 2},\\
& \vert U\vert^{\otimes m}z_m + \text{terms in $z_2, z_3, \dots,
z_{m-1}$} = 0
\,\in H^{\otimes m}, \quad \text{for $m \geq 3$.}
\endsplit
$$
We deduce $z_m = 0$ by induction on degree $m$, 
since $\vert U\vert^{\otimes m}$ is an isomorphism for each $m$.
Hence $U: \T \to \T$ is an $R$-algebra isomorphism. 
If $w$ in the equation (1.4) satisfies $w_2 = w_3 = \dots = 
w_{p-1} = 0$, then the solution $z$ also satisfies 
$z_2 = z_3 = \dots = z_{p-1} = 0$ 
because $\vert U\vert^{\otimes m}$ is an isomorphism. 
This means $U\inv(\T_p) \subset \T_p$ for each $p \geq 2$. 
It is also proved $U\inv(\T_1) \subset \T_1$ in a similar way. 
Consequently $U$ is a filtration-preserving $R$-algebra 
automorphism of the algebra $\T$, that is, $U \in \Aut(\T)$.\par
This completes the proof of the lemma.\qed
\enddemo

As is the lemma, we denote by $\vert U\vert \in \GL(H)$ 
the automorphism of $H = \T_1/\T_2$ induced by the 
automorphism $U \in \Aut(\T)$. This defines a group 
homomorphism 
$$
\vert\cdot\vert: \Aut(\T) \to \GL(H), \quad
U \mapsto \vert U\vert,
$$ 
whose kernel we denote by 
$$
\IA(\T) := \Ker \vert\cdot\vert \subset \Aut(\T). \tag 1.5
$$
Any element $A \in \GL(H)$ can be regarded as 
a filtration-preserving automorphism of $\T$ by 
$A(z_m) = (A^{\otimes m}z_m)$, $z_m \in H^{\otimes m}$.
Hence we have a semi-direct product decomposition 
$$
\Aut(\T) = \IA(\T) \rtimes \GL(H). \tag 1.6
$$
From Lemma 1.2 we have a natural bijection 
$$
E: \IA(\T) \to \Hom(H, \T_2), \quad
U \mapsto  U_H - 1_H.
\tag 1.7
$$
We often identify $\IA(\T)$ and $\Hom(H, \T_2)$ by the map $E$. 
\par

The group $\IA(\T)$ acts on the set $\Theta_n$ 
of all the Magnus expansions in a natural way. 
If $U \in \IA(\T)$ and $\theta \in \Theta_n$, 
then the composite $U \circ \theta: \Fn\overset\theta\to\to
1+\T_1\overset{U}\to\to1+\T_1$ is also a Magnus expansion. 
To study the action, recall the completed group ring $\RFn$ 
of the group $\Fn$. By definition, it is the completion 
of the group ring $R[\Fn]$ 
with respect to the augmentation ideal $I_R[\Fn] =
\{\sum_{\gamma \in \Fn}r_\gamma\gamma \in R[\Fn]; \, 
\sum r_\gamma = 0\}$
$$
\RFn := \varprojlim_{m\to +\infty} R[\Fn]/I_R[\Fn]^m,
$$
which has a decreasing filtration of the ideals
$$
\widehat{I_R[\Fn]^p\,} := \varprojlim_{p \leq m\to +\infty}
I_R[\Fn]^p/I_R[\Fn]^m, \quad p \geq 1.
$$
Any Magnus expansion $\theta \in \Theta_n$ induces an
$R$-algebra homomorphism $\theta: R[\Fn] \to \T$ in an obvious
way. Since $\theta(I_R[\Fn]) \subset \T_1$, 
we obtain an $R$-algebra homomorphism $\theta: \RFn \to \T$, 
which maps the ideal $\widehat{I_R[\Fn]^p\,}$ into $\T_p$ for
each $p\geq 1$.\par
Our construction of Johnson maps of the automorphism group 
$\Aut(\Fn)$ is based on 
\proclaim{Theorem 1.3}{\rm (1)} For any $\theta \in \Theta_n$ 
the homomorphism
$$
\theta: \RFn \to \T
$$
is an $R$-algebra isomorphism, which maps $\widehat{I_R[\Fn]^p\,}$
\underbar{onto} $\T_p$ for each $p \geq 1$.\par
{\rm (2)} If $\theta'$ and $\theta'' \in \Theta_n$, then 
there exists a unique $U \in \IA(\T)$ such that
$$
\theta'' = U\circ\theta' \in \Theta_n.
$$
In other words, the action of the group $\IA(\T)$ on the set 
$\Theta_n$ is free and transitive.
\endproclaim

\demo{Proof} First we prove the assertion (1) for the standard 
Magnus expansion $\std: \Fn \to 1+\T_1$, $\std(x_i) = 1+X_i$, 
$1 \leq i \leq n$. The $R$-algebra homomorphism $\kappa: \T
\to \RFn$ given by $\kappa(X_i) = x_i - 1$ maps $\T_p$ into 
$\widehat{I_R[\Fn]^p\,}$ for each $p$. Now we have 
$$
\split
& \std\circ\kappa(X_i) = \std(x_i-1) = 1+X_i-1 = X_i, 
\quad\text{and}\\
& \kappa\circ\std(x_i) = \kappa(X_i) = 1+x_i-1 = x_i.
\endsplit
$$
These imply $\std\circ\kappa = 1_\T$ and $\kappa\circ\std = 
1_{\RFn}$, respectively. 
Hence $\std: \RFn \to \T$ is an $R$-algebra isomorphism, 
which maps $\widehat{I_R[\Fn]^p\,}$ \underbar{onto} $\T_p$ 
for each $p$.\par
Next we consider an arbitrary Magnus expansion 
$\theta \in \Theta_n$. The $R$-algebra endomorphism 
$\theta\circ\kappa:
\T\overset\kappa\to\to\RFn\overset\theta\to\to \T$ 
satisfies the conditions in (i) and (ii) in Lemma 1.2. 
In fact, $\vert\theta\circ\kappa\vert = 1_H$. 
Hence $\theta\circ\kappa \in \IA(\T)$, which we denote by $U$.
Especially $\theta = U \circ\std: 
\RFn \to \T$ is an $R$-algebra isomorphism, which maps
$\widehat{I_R[\Fn]^p\,}$ \underbar{onto} $\T_p$. 
This implies the action of $\IA(\T)$ on $\Theta_n$ is transitive.
\par
Finally we prove the action is free. Suppose $U \in \IA(\T)$ 
satisfies $U \circ\std = \std: \Fn \to 1 + \T_1$. 
Then we have $U \circ\std = \std: \RFn \overset\cong\to\to \T$, 
and so $U = U \circ\std\circ\kappa = \std\circ\kappa = 1_\T$, 
as was to be shown.\par
This completes the proof of the theorem.\qed
\enddemo

We conclude the section by writing down the group structure 
on the set $\Hom(H, \allowmathbreak \T_2)\times \GL(H)$ induced 
by the decomposition (1.6) and the bijection (1.7) in low degree.  
We denote $((u, A)) := (E\inv u)\circ A \in \Aut(\T)$
for $u \in \Hom(H, \T_2)$ and $A \in \GL(H)$. We have
$$
((u, A))a = Aa + {\sum}^\infty_{m=1}u_m(Aa)
$$
for any $a \in H$.
By straightforward computation we obtain
\proclaim{Lemma 1.4} Suppose
$((w, C)) = ((u, A))((v, B)) \in \Aut(\T)$
for $(u, A)$, $(v, B)$ and $(w, C) \in 
\Hom(H, \T_2)\times \GL(H)$.
Then we have 
$$
\align
& C = AB,\\
& w_1 = u_1 + Av_1, \quad\text{and}\\
& w_2 = u_2 + (u_1\otimes 1 + 1\otimes u_1)Av_1 + Av_2,
\endalign
$$
where $u_p$, $v_p$ and $w_p \in \Hom(H, H^{\otimes (p+1)})$ 
are the $p$-th components of $u$, $v$ and $w$, respectively, and
$$
Av_p := (\oversetbrace{p+1}\to{A\otimes \cdots\otimes A})v_pA\inv
\in \Hom(H, H^{\otimes(p+1)}). 
$$
\endproclaim
\demo{Proof} For any $a \in H$ we have 
$$
\split
& Ca + w_1(Ca) + w_2(Ca) 
\equiv ((w, C))a = ((u, A))((v, B))a\\
\equiv\, & ((u, A))(Ba + v_1(Ba) + v_2(Ba))\\
\equiv\, & ABa + u_1(ABa) + u_2(ABa)\\
& \quad + (A\otimes A)v_1(Ba) + (u_1\otimes 1 + 1\otimes
u_1)(A\otimes A)v_1(Ba)\\
& \quad + (A\otimes A\otimes A)v_2(Ba)\\
=\, & ABa + (u_1+Av_1)(ABa) + (u_2 + (u_1\otimes 1 + 1\otimes
u_1)Av_1 + Av_2)(ABa)
\endsplit
$$
modulo $\T_4$, as was to be shown.\qed
\enddemo

\beginsection 2. Johnson Maps.\par

As was shown in Theorem 1.3(2), the group $\IA(\T)$ acts on the 
set $\Theta_n$ of all the $R$-valued Magnus expansions of the 
free group $\Fn$ in a free and transitive way. 
We denote by $H^*$ the dual of $H$, $H^* := \Hom(H, R)$. 
We often identify 
$$
\IA(\T) \overset{E}\to\cong \Hom(H, \T_2) 
= {\prod}^\infty_{p=1} \Hom(H, H^{\otimes (p+1)})
= {\prod}^\infty_{p=1}H^*\otimes H^{\otimes (p+1)}.
\tag 2.1
$$
Now we consider the automorphism group of the group $\Fn$, 
$\Aut(\Fn)$. It acts on the set $\Theta_n$ in a natural way. 
In fact, we define 
$$
\varphi\cdot\theta :=
\vert\varphi\vert\circ\theta\circ\varphi\inv
$$
for $\varphi \in \Aut(\Fn)$ and $\theta \in \Theta_n$. 
Here $\vert\varphi\vert \in \GL(H)$ is the induced map 
on $H = H_1(\Fn; R)$ by the automorphism $\varphi$. 
From the free and transitive action of $\IA(\T)$, 
there exists a unique automorphism $\tau^\theta(\varphi) 
\in \IA(\T)$ such that 
$$
\varphi\cdot\theta = {\tau^\theta(\varphi)}\inv\circ\theta 
\in \Theta_n. \tag 2.2
$$
If we fix a Magnus expansion $\theta \in \Theta_n$, 
it defines a map
$$
\tau^\theta: \Aut(\Fn) \to \IA(\T), \quad 
\varphi\mapsto \tau^\theta(\varphi),
$$
which we call {\it the total Johnson map induced by the 
Magnus expansion $\theta$}. Immediately from (2.2) we have 
a commutative diagram
$$
\CD
\RFn @>{\theta}>> \T\\
@V{\varphi}VV @V{\tau^\theta(\varphi)\circ\vert\varphi\vert}VV\\
\RFn @>{\theta}>> \T.\\
\endCD \tag 2.3
$$
Hence we obtain 
$$
\tau^\theta(\varphi\psi) =
\tau^\theta(\varphi)\circ\vert\varphi\vert\circ\tau^\theta(\psi)
\circ\vert\varphi\vert\inv \tag 2.4
$$
for any $\varphi$ and $\psi \in \Aut(\Fn)$. \par

Under the identification (2.1), for each $p \geq 1$, we define {\it 
the $p$-th
Johnson map  induced by the $R$-valued Magnus expansion $\theta$}
$$
\tau^\theta_p: \Aut(\Fn) \to H^*\otimes H^{\otimes (p+1)}, \quad
\varphi \mapsto \tau^\theta_p(\varphi)
$$
by the $p$-th component 
of the total Johnson map $E\tau^\theta$. We have 
$$
E\tau^\theta(\varphi) = {\sum}^\infty_{p=1}\tau^\theta_p(\varphi)
\in \Hom(H, \T_2) = {\prod}^\infty_{p=1}H^*\otimes H^{\otimes (p+1)}
$$
for $\varphi \in \Aut(\Fn)$.\par
The map $\tau^\theta_p$ is \underbar{not} a group 
homomorphism. The relation (2.4) means an infinite sequence of 
coboundary relations. 
In this paper we confine ourselves to studying the first and 
the second relations, which have cohomological consequences 
about the group $\An$ and the mapping class groups for surfaces. 
In the case $p = 1$ and $2$, from Lemma 1.4, we have
$$
\align
& \tau^\theta_1(\varphi\psi) 
= \tau^\theta_1(\varphi) + \vert\varphi\vert\tau^\theta_1(\psi)\\
& \tau^\theta_2(\varphi\psi) 
= \tau^\theta_2(\varphi) 
+ (\tau^\theta_1(\varphi)\otimes 1 +
1\otimes\tau^\theta_1(\varphi))
\vert\varphi\vert\tau^\theta_1(\psi)
+ \vert\varphi\vert\tau^\theta_2(\psi)\\
\endalign
$$ 
for any $\varphi$ and $\psi \in \An$.
In the succeeding sections we will show 
these elementary formulae have some significant consequences 
in the cohomology of the group $\An$.
Throughtout this paper we denote by $C^*(G; M)$ 
the normalized standard complex of a group $G$ with values in 
a $G$-module $M$, and use the Alexander-Whitney cup product 
$\cup: C^*(G; M_1) \otimes C^*(G; M_2)\to C^*(G; M_1\otimes
M_2)$. For details, see [HS] ch.II. The formulae are equivalent to 

\proclaim{Lemma 2.1}
$$
\align
& -d\tau^\theta_1 = 0 \,\in\, C^*(\An; H^*\otimes H^{\otimes 2}),
\tag 2.5\\
& -d\tau^\theta_2 = 
(\tau^\theta_1\otimes 1 + 1\otimes\tau^\theta_1) \cup 
\tau^\theta_1 \,\in\, C^*(\An; H^*\otimes H^{\otimes 3}).
\tag 2.6\\
\endalign
$$
In (2.6) we drop the composite map
$\Hom(H^{\otimes 2}, H^{\otimes 3})\otimes \Hom(H, H^{\otimes 2})
\to \Hom(H, \allowmathbreak H^{\otimes 3}) = H^*\otimes H^{\otimes 3}$, 
$f\otimes g \mapsto f\circ g$, for simplicity.
\endproclaim

From (2.5) the map $\tau^\theta_1$ is a $1$-cocycle 
of the group $\An$ with values in the $\An$-module 
$H^*\otimes H^{\otimes 2}$. The cohomology class 
$[\tau^\theta_1] \in H^1(\An; H^*\otimes H^{\otimes 2})$
is independent of the choice of a Magnus expansion $\theta$. 
It can be proved directly from Theorem 1.3(2). 
As will be shown in \S4, the class $[\tau^\theta_1]$ is 
the Gysin image of a certain cohomology class in 
$H^2(\Fn \rtimes \An; H^{\otimes 2})$ independent of 
the choice of a Magnus expansion. \par

\proclaim{Lemma 2.2} We have 
$$
\align
&\tau^\theta_1(\varphi)\vert\varphi\vert[\gamma] 
= \theta_2(\varphi(\gamma)) 
-\vert\varphi\vert^{\otimes 2}\theta_2(\gamma) \tag 2.7\\
&\tau^\theta_1(\varphi)[\gamma] 
= \theta_2(\gamma) 
-\vert\varphi\vert^{\otimes 2}\theta_2(\varphi\inv(\gamma)) \tag 2.8\\
\endalign
$$
for any $\gamma \in \Fn$ and $\varphi \in \An$.
\endproclaim
In fact, from (2.3), we have
$$
\split
& 1 + [\varphi(\gamma)] + \theta_2(\varphi(\gamma)) 
\equiv \theta(\varphi(\gamma)) 
= \tau^\theta(\varphi)\vert\varphi\vert\theta(\gamma)\\
\equiv\,& \tau^\theta(\varphi)(1+\vert\varphi\vert[\gamma] 
+ \vert\varphi\vert^{\otimes 2}\theta_2(\gamma))\\
\equiv\,& 1+\vert\varphi\vert[\gamma] 
+ \tau^\theta_1(\varphi)\vert\varphi\vert[\gamma]
+ \vert\varphi\vert^{\otimes 2}\theta_2(\gamma)
\endsplit
$$
modulo $\T_3$.\par

Finally we compute the Johnson maps 
on the inner automorphisms of the group $\Fn$. 
The image of the homomorphism
$$
\iota: \Fn \to \An, \quad \gamma\mapsto
\left(\iota(\gamma): \delta\mapsto \gamma\delta\gamma\inv\right)
\tag 2.9
$$
is, by definition, the inner automorphism group of $\Fn$, 
and often denoted by $\Inn$. 
The quotient $\On := \An/\Inn$ is called 
the outer automorphism group of $\Fn$.\par

\proclaim{Lemma 2.3} For $p \geq 1$, $\gamma \in \Fn$ and 
$a \in H$ we have 
$$
\tau^\theta_p(\iota(\gamma))a = \theta_p(\gamma)a + 
\sum^m_{j=1}\sum_{q_0 + q_1 +\cdots q_j = p\atop 
q_0 \geq 0,\,\, q_1, \dots, q_j \geq 1}(-1)^j
\theta_{q_0}(\gamma)\,a\,\theta_{q_1}(\gamma)\cdots 
\theta_{q_j}(\gamma).
\tag 2.10
$$
\endproclaim
\demo{Proof} Recall
$$
\theta(\gamma)\inv
= \left(1 +{\sum}^\infty_{s=1}\theta_s(\gamma)\right)\inv
= 1 +
{\sum}^\infty_{j=1}(-1)^j
\left({\sum}^\infty_{s=1}\theta_s(\gamma)\right)^j.
$$
The map $U_\gamma: \T \to \T$, 
$z \mapsto \theta(\gamma)z\theta(\gamma)\inv$, 
is an element of $\Aut(\T)$. 
Now, since $\vert\iota(\gamma)\vert = 1$, we have
$$
\split
& \tau^\theta(\iota(\gamma))a = U_\gamma a 
= \theta(\gamma)a\theta(\gamma)\inv\\
& = \left(\sum^\infty_{q_0=0}\theta_{q_0}(\gamma)\right)a
\left(1 +
\sum^\infty_{j=1}(-1)^j\sum^\infty_{q_1,
\dots,
q_j=1}\theta_{q_1}(\gamma)\cdots \theta_{q_j}(\gamma)\right).
\endsplit
$$
Taking the $(p+1)$-st components in $H^{\otimes (p+1)}$, 
we obtain the lemma.\qed
\enddemo

In the case $p=1$ and $2$, we have
$$
\align
& \tau^\theta_1(\iota(\gamma))a = [\gamma]a-a[\gamma],\tag 2.11\\
& \tau^\theta_2(\iota(\gamma))a =
\theta_2(\gamma)a - a\theta_2(\gamma) +
a[\gamma][\gamma]-[\gamma]a[\gamma].
\tag 2.12\\
\endalign
$$

\beginsection 3. Lower Central Series. \par

In this section we suppose that the natural ring homomorphism 
$\nu \in \bZ \mapsto \nu\cdot 1\in R$ is injective. 
We denote by $\Gamma_m = \Gamma_m(\Fn)$, $m \geq 1$, 
the lower central series of the group $\Fn$
$$
\Gamma_1 := \Fn, \qquad
\Gamma_{m+1} := [\Gamma_m, \Fn], \quad m \geq 1.
$$
The group $\Aut(\Fn)$ acts on the subgroup $\Gamma_m$ and 
the quotient $\Gamma_1/\Gamma_m$ in a natural way.
S\. Andreadakis [An] introduced a decreasing filtration $\{A(m)\}^\infty_{m=1}$
of $\Aut(\Fn)$ by 
$$
A(m) := \Ker(\Aut(\Fn) \to \Aut(\Gamma_1/\Gamma_{m+1})), \quad
m \geq 0.
$$
In [An] he wrote $K_m$ for $A(m)$. 
The $m$-th Johnson homomorphism $\tau_m = \tau^J_m$ describes 
the quotient $A(m)/A(m+1)$, which was introduced by 
D\. Johnson [J1]. The homomorphism $\tau_m$ can be regarded as 
an embedding of the quotient $A(m)/A(m+1)$ into the module 
$H^*\otimes H^{\otimes (m+1)}$. We prove the restriction of 
$\tau^\theta_m$ to $A(m)$ coincides with $\tau_m$ (Theorem 3.1). 
Especially $\tau^\theta_m\vert_{A(m)}$ is a homomorphism independent 
of the choice of a Magnus expansion $\theta$. \par

Choose a Magnus expansion $\theta \in \Theta_n$. As was proved by Magnus 
[M3], we have
$$
\theta\inv(1+\T_m) =\Gamma_m \tag 3.1
$$
for each $m \geq 1$. See [Bou] ch\. 2, \S5, no\. 4, Theorem 2. 
This implies the $m$-th component $\theta_m$ gives an injective homomorphism
$$
\theta_m: \Gamma_m/\Gamma_{m+1} \hookrightarrow H^{\otimes m}.
$$
\par
Here the restriction $\theta_m\vert_{\Gamma_m}$ is independent of 
the choice of $\theta$. We prove it by induction on $m$. 
From Definition 1.1, $\theta_1$ is independent of $\theta$. 
Assume $m \geq 2$. We have $\Gamma_m = [\Gamma_{m-1}, \Gamma_1]$. 
Let $\gamma \in \Gamma_{m-1}$. From (3.1) follows 
$\theta(\gamma) \equiv 1 + \theta_{m-1}(\gamma) \pmod{\T_{m}}$. 
For $\delta \in \Gamma_1$, we have, modulo $\T_{m+1}$, 
$$
\split
& \theta(\gamma)\theta(\delta)\theta(\gamma\inv)\theta(\delta\inv)\\
= \,& \theta(\gamma)(1+(\theta(\delta)-1))\theta(\gamma\inv)\theta(\delta\inv)\\
= \,&
\theta(\delta\inv)
+ \theta(\gamma)(\theta(\delta)-1)\theta(\gamma\inv)\theta(\delta\inv)\\
\equiv \,&
\theta(\delta\inv) + (1+\theta_{m-1}(\gamma))
(\theta(\delta)-1)(1-\theta_{m-1}(\gamma))\theta(\delta\inv)\\
\equiv \,& \theta(\delta\inv) 
+ (\theta(\delta)-1)\theta(\delta\inv) 
+ \theta_{m-1}(\gamma)(\theta(\delta)-1)\theta(\delta\inv)
- (\theta(\delta)-1)\theta_{m-1}(\gamma)\theta(\delta\inv)\\
\equiv \,& 1+ \theta_{m-1}(\gamma)[\delta] - [\delta]\theta_{m-1}(\gamma).\\
\endsplit
$$ 
Hence
$$
\theta_m(\gamma\delta\gamma\inv\delta\inv) = 
\theta_{m-1}(\gamma)[\delta] -[\delta]\theta_{m-1}(\gamma). \tag 3.2
$$
From the inductive assumption $\theta_{m-1}(\gamma)$ is independent of 
the choice of $\theta$. This completes the induction.\par

We denote the image $\theta_m(\Gamma_m/\Gamma_{m+1})$ 
by $\Cal L_m \subset H^{\otimes m}$. 
We identify $\Gamma_m/\Gamma_{m+1}$ 
and $\Cal L_{m}$ by the isomorphism $\theta_m$. As is known, 
the sum $\bigoplus^\infty_{m=1} \Cal L_m \subset \T$ is a Lie subalgebra 
of the associative algebra $\T$, and naturally isomorphic 
to the free Lie algebra $\Cal L(H_\bZ)$ generated by $H_\bZ = {\Fn}\abel$. 
See [Bou] loc\. cit. \par
Now we recall the definition of the Johnson homomorphisms [J1][J2]. 
If $\varphi \in A(m)$ and $\gamma \in \Fn$, then we have 
$\gamma\inv\varphi(\gamma) \in \Gamma_{m+1}$. This allows us to consider 
$$
\tau_m(\varphi, \gamma) := \gamma\inv\varphi(\gamma) \mod \Gamma_{m+2} 
\in \Cal L_{m+1}.
$$
It is easy to prove the map $\gamma \in \Fn \mapsto 
\tau_m(\varphi, \gamma) \in \Cal L_{m+1}$ is a group homomorphism. 
Hence it can be regarded as an element $\tau_m(\varphi) 
\in \Hom_\bZ(H_\bZ, \Cal L_{m+1})$, 
and induces a map
$$
\tau_m = \tau^J_m: A(m) \to \Hom_\bZ(H_\bZ, \Cal L_{m+1}), \quad
\varphi \mapsto \tau_m(\varphi).
$$
Moreover one can easily prove $\tau_m$ is a group homomorphism. 
The homomorphism $\tau_m$ is, by definition, 
the $m$-th Johnson homomorphism [J2]. 
Immediately from the definition we have
$$
\Ker\, \tau_m = A(m+1), \tag 3.3
$$
so that it can be regarded as an embedding 
$$
\tau_m: A(m)/A(m+1) 
\hookrightarrow \Hom_\bZ(H_\bZ, \Cal L_{m+1}).
$$
\par
Since the integers $\bZ$ is a subring of $R$, we may regard 
$\Hom_\bZ(H_\bZ, \Cal L_{m+1})$ as a $\bZ$-submodule of 
$H^*\otimes H^{\otimes(m+1)}$ in an obvious way. Then 
\proclaim{Theorem 3.1} We have 
$$
\tau_m = \tau^\theta_m\vert_{A(m)}: A(m) \to H^*\otimes H^{\otimes (m+1)}
$$
for each  $m \geq 1$. 
Especially the restriction $\tau^\theta_m\vert_{A(m)}$ 
is a group homomorphism independent of the choice of 
the Magnus expansion $\theta$. 
\endproclaim
\demo{Proof} We prove the theorem by induction on $m \geq 1$. 
For any $\varphi \in A(m) \subset A(1)$ we have $\vert\varphi\vert = 1$. 
If $m=1$, then we have, modulo $\T_3$, $\theta(\varphi(\gamma)) =
\tau^\theta(\varphi)\theta(\gamma) \equiv \theta(\gamma) +
\tau_1^\theta(\varphi)[\gamma]$, and so
$$
\split
& 1 + \tau_1(\varphi, \gamma) = 1 + \theta_2(\gamma\inv\varphi(\gamma))
\equiv \theta(\gamma\inv\varphi(\gamma)) =
\theta(\gamma)\inv\theta(\varphi(\gamma))\\
 \equiv\,&
\theta(\gamma)\inv(\theta(\gamma) +
\tau_1^\theta(\varphi)[\gamma]) \equiv 1 + \tau_1^\theta(\varphi)[\gamma].
\endsplit
$$
This implies we have $\tau^\theta_1\vert_{A(2)} = 0$ from (3.3).\par
Suppose $m \geq 2$. From the inductive assumption and (3.3) we have 
$\tau^\theta_1(\varphi) = \cdots = \tau^\theta_{m-1}(\varphi) = 0$ 
for any $\varphi \in A(m)$, and so $\theta(\varphi(\gamma)) \equiv 
\theta(\gamma) + \tau^\theta_m(\varphi)[\gamma] \pmod{\T_{m+2}}$. 
Hence we have 
$$
\split
& 1 + \tau_m(\varphi, \gamma) = 1 + \theta_{m+1}(\gamma\inv\varphi(\gamma))
\equiv \theta(\gamma\inv\varphi(\gamma)) =
\theta(\gamma)\inv\theta(\varphi(\gamma))\\
 \equiv\,&
\theta(\gamma)\inv(\theta(\gamma) +
\tau_m^\theta(\varphi)[\gamma]) \equiv 1 + \tau_m^\theta(\varphi)[\gamma]
\endsplit
$$
modulo $\T_{m+2}$. This completes the induction and 
the proof of the theorem.\qed
\enddemo

\beginsection 4. Twisted Cohomology Classes. \par

In this section we introduce two series of twisted cohomology classes 
$$
h_p \in H^p(\An; H^*\otimes H^{\otimes (p+1)}) \quad \text{and}\quad
\h_p \in H^p(\An; H^{\otimes p}) 
$$
for $p \geq 1$ by an analogous construction 
to the Morita-Mumford classes 
on the mapping class groups for surfaces [Mu] [Mo1]. 
Restricted to the mapping class group $\M1g$ of genus $g$ 
with $1$ boundary component, they coincide with 
the twisted Morita-Mumford classes [Ka] [KM1] [KM2]
$$
\align
& (p+2)!\,h_p\vert_{\M1g} = m_{0, p+2} \in H^p(\M1g; H^{\otimes (p+2)}),
\quad\text{and}\tag 5.8\\
& p!\,\h_p\vert_{\M1g} = -m_{1, p} \in H^p(\M1g; H^{\otimes p}),\tag 5.18
\endalign
$$ 
as will be shown in \S5. 
We prove a suitable algebraic combination of $p$ copies of 
the $1$-cocycle $\tau^\theta_1$ introduced in \S2  
represents the cohomology class $h_p$ for each $p \geq 1$
(Theorem 4.1). In the case $p=1$, $\tau^\theta_1$ represents the class $h_1$.
Furthermore we describe some contraction formulae deduced 
from the relation (2.6).\par

In order to define the cohomology classes $h_p$ and $\h_p$, 
we consider the semi-direct product 
$$
\Anbar := \Fn\rtimes\An
$$ 
and the map
$$
k_0: \Anbar \to H, \quad (\gamma, \varphi) \mapsto [\gamma], \tag4.1
$$
introduced in [Mo3]. 
The group $\Anbar$ is, by definition, the product set $\Fn\times\An$ 
with the group law
$$
(\gamma_1, \varphi_1) (\gamma_2, \varphi_2) 
:= (\gamma_1\varphi_1(\gamma_2), \varphi_1\varphi_2), \quad
(\gamma_i \in \Fn, \, \varphi_i \in \An).
$$
We often write simply $\gamma\varphi$ for $(\gamma, \varphi)$. 
It is easy to prove $k_0$ satisfies the cocycle condition. 
We write also $k_0$ for the cohomology class 
$[k_0] \in H^1(\Anbar; H)$. 
Consider the $(p+1)$-st power of $k_0$
$$
{k_0}^{\otimes (p+1)} \in H^{p+1}(\Anbar; H^{\otimes (p+1)})
$$
for each $p \geq 0$.\par
The group $\Anbar$ admits a group extension 
$$
\Fn \overset\iota\to\to \Anbar \overset\pi\to\to \An \tag 4.2
$$
given by $\iota(\gamma) = \gamma$ and $\pi(\gamma\varphi) = \varphi$ 
for $\gamma \in \Fn$ and $\varphi \in \An$. It induces the Gysin map 
$$
\pi_\sharp: H^{p+1}(\Anbar; H^{\otimes (p+1)})
\to H^{p}(\An; H^{*}\otimes H^{\otimes (p+1)}).
$$
Here we identify
$$
H^1(\Fn; H^{\otimes (p+1)}) = \Hom(H, H^{\otimes (p+1)})
= H^{*}\otimes H^{\otimes (p+1)} \tag 4.3
$$
in a natural way. 
For each $p \geq 1$ we define
$$
h_p := \pi_\sharp({k_0}^{\otimes (p+1)}) \in H^p(\An; 
H^*\otimes H^{\otimes (p+1)}). \tag 4.4
$$
In the case $p=0$ we have 
$$
\pi_\sharp(k_0) = \iota^*k_0 = 1_H \in H^0(\An; H^*\otimes H).
\tag 4.5
$$
Contracting the coefficients by the $\GL(H)$-homomorphism
$$
\a_p: H^*\otimes H^{\otimes (p+1)} 
\to H^{\otimes p},\quad
f\otimes v_0\otimes v_1 \otimes \cdots \otimes v_p \mapsto 
f(v_0) v_1 \otimes \cdots \otimes v_p,
\tag 4.6
$$
we define
$$
\h_p := {\a_p}_*(h_p) \in H^p(\An; H^{\otimes p}).\tag 4.7
$$
\par

We introduce a $\GL(H)$-homomorphism
$$
\bb_p: (H^*\otimes H^{\otimes 2})^{\otimes p} = 
\Hom(H, H^{\otimes 2})^{\otimes p} \to \Hom(H, H^{\otimes (p+1)})
= H^*\otimes H^{\otimes (p+1)}
$$
for each $p \geq 1$. If $p \geq 2$, we define
$$
\aligned
& \bb_p(u_{(1)}\otimes u_{(2)}\otimes
\cdots\otimes u_{(p-1)}\otimes u_{(p)}) \\
& := 
\left(u_{(1)}\otimes {1_H}^{\otimes (p-1)}\right)\circ
\left(u_{(2)}\otimes {1_H}^{\otimes (p-2)}\right)\circ\cdots\circ
\left(u_{(p-1)}\otimes {1_H}\right)\circ u_{(p)},
\endaligned
\tag4.8
$$
where $u_{(i)} \in \Hom(H, H^{\otimes 2}) = H^*\otimes H^{\otimes 2}$, 
$1 \leq i \leq p$.  
In the case $p=1$, we define $\bb_1 := 1_{H^*\otimes H^{\otimes 2}}$.
Now we prove
\proclaim{Theorem 4.1} 
$$
h_p = {\bb_p}_*([\tau^\theta_1]^{\otimes p}) 
\in H^p(\An; H^*\otimes H^{\otimes(p+1)})
$$
for any Magnus expansion $\theta \in \Theta_n$ and each $p \geq 1$. 
In the case $p=1$ we have $[\tau^\theta_1] = h_1
\in  H^1(\An; H^*\otimes H^{\otimes 2})$, which is independent of the 
choice of $\theta$.
\endproclaim
\demo{Proof}
We define a $1$-cochain $\widetilde{\theta}_2 \in C^1(\Anbar; H^{\otimes 2})$ by 
$\widetilde{\theta}_2(\gamma\varphi) := \theta_2(\gamma)$ 
for $\gamma \in \Fn$ and $\varphi \in \An$. 
Then we have
$$
d\widetilde{\theta}_2 = -\left(\tau^\theta_1\circ k_0 + {k_0}^{\otimes 2}\right) 
\in C^2(\Anbar; H^{\otimes 2}). \tag4.9
$$
In fact, it follows from (2.7) 
$$
\split
& d\widetilde{\theta}_2(\gamma_1\varphi_1, \gamma_2\varphi_2)
=  \vert\varphi_1\vert^{\otimes 2}\widetilde{\theta}_2(\gamma_2\varphi_2) 
- \widetilde{\theta}_2(\gamma_1\varphi_1(\gamma_2)\varphi_1\varphi_2) 
+ \widetilde{\theta}_2(\gamma_1\varphi_1)\\
= \, & \vert\varphi_1\vert^{\otimes 2}\theta_2(\gamma_2) 
- \theta_2(\gamma_1\varphi_1(\gamma_2)) 
+ \theta_2(\gamma_1)\\
= \, & \vert\varphi_1\vert^{\otimes 2}\theta_2(\gamma_2) 
- \theta_2(\varphi_1(\gamma_2)) 
- [\gamma_1]\otimes [\varphi_1(\gamma_2)]\\
= \, & -\tau^\theta_1(\varphi_1)\vert\varphi_1\vert k_0(\gamma_2\varphi_2)
- k_0(\gamma_1\varphi_1)\otimes \vert\varphi_1\vert
k_0(\gamma_2\varphi_2)\\ 
= \, & -\left(\tau^\theta_1\circ k_0 + {k_0}^{\otimes 2}\right)
 (\gamma_1\varphi_1, \gamma_2\varphi_2).
\endsplit
$$
Consider a $(p+1)$-cocycle $f_p \in C^{p+1}(\Anbar; H^{\otimes (p+1)})$ defined
by 
$$
\aligned
f_p & := \left({\bb_p}_*(\tau^\theta_1)^{\otimes p}\right)\circ k_0\\
& = \left(\tau^\theta_1\otimes{1_H}^{\otimes(p-1)}\right)\circ
\left(\tau^\theta_1\otimes{1_H}^{\otimes(p-2)}\right)\circ \cdots \circ
\left(\tau^\theta_1\otimes{1_H}\right)\circ
\tau^\theta_1\circ k_0
\endaligned
$$
and a $p$-cochain $g_p \in C^{p}(\Anbar; H^{\otimes (p+1)})$ defined by 
$$
g_p := 
\left(\tau^\theta_1\otimes{1_H}^{\otimes(p-1)}\right)\circ
\left(\tau^\theta_1\otimes{1_H}^{\otimes(p-2)}\right)\circ \cdots \circ
\left(\tau^\theta_1\otimes{1_H}\right)\circ\widetilde{\theta}_2
$$
for $p \geq 1$. If $p = 0$, we define $f_0 := k_0$. 
From (4.9) follows
$$
\split
d g_p &= (-1)^{p-1}\left(\tau^\theta_1\otimes{1_H}^{\otimes(p-1)}\right)
\circ \cdots \circ
\left(\tau^\theta_1\otimes{1_H}\right)\circ d\widetilde{\theta}_2\\
&= (-1)^{p}\left(\tau^\theta_1\otimes{1_H}^{\otimes(p-1)}\right)\circ
\cdots \circ
\left(\tau^\theta_1\otimes{1_H}\right)\circ
\left(\tau^\theta_1\circ k_0 + {k_0}^{\otimes 2}\right)\\
&= (-1)^{p}\left(f_p + f_{p-1}\otimes k_0\right).
\endsplit
$$
Hence we obtain
$$
\split
&[\left({\bb_p}_*(\tau^\theta_1)^{\otimes p}\right)\circ k_0] = [f_p] =
-[f_{p-1}\otimes k_0] = \cdots  \\
= \,&(-1)^{p-1}[f_1\otimes {k_0}^{\otimes
(p-1)}]  = (-1)^{p}[{k_0}^{\otimes (p+1)}],
\endsplit
$$
that is,
$$
\left({\bb_p}_*(\tau^\theta_1)^{\otimes p}\right)\circ k_0 =
(-1)^{p}{k_0}^{\otimes (p+1)} \in H^{p+1}(\Anbar; H^{\otimes (p+1)}),
\tag4.10
$$
for each $p \geq 1$. We have 
$$
\aligned
h_p =\, &\pi_\sharp({k_0}^{\otimes (p+1)}) = 
{\bb_p}_*(\tau^\theta_1)^{\otimes p}\circ \pi_\sharp k_0
= {\bb_p}_*(\tau^\theta_1)^{\otimes p}
\in H^p(\An; H^*\otimes H^{\otimes (p+1)}).
\endaligned\tag 4.11
$$
In fact, the cocycle $f_{p} \in C^{p+1}(\Anbar; H^{\otimes (p+1)})$ is
contained  in the $p$-th filter $A_p$ introduced in [HS] ch.II, p.118.
Therefore,  from [HS] ch.II, Proposition 3, p.125, its Gysin image is given by
$$
\split
&\left(\pi_\sharp f_{p}\right)(\varphi_1, \dots, \varphi_p)\,\,[\gamma] 
= \left(f_{p}\right)_1(\gamma, \varphi_1, \dots, \varphi_p) \\
= & (-1)^pf_{p}(\varphi_1, \dots, \varphi_p, 
\left(\varphi_1\cdots\varphi_p\right)\inv(\gamma))\\
= &(-1)^p\left(\left(\left(\tau^\theta_1\otimes{1_H}^{\otimes(p-1)}\right)\circ
\cdots \circ\left(\tau^\theta_1\otimes{1_H}\right)\circ
\tau^\theta_1\right)(\varphi_1, \dots, \varphi_p)\right)\circ 
\vert\varphi_1\cdots\varphi_p \vert\\
& \phantom{(-1)^p\left(\left(\tau^\theta_1\otimes{1_H}^{\otimes(p-1)}\right)\circ
\cdots \circ\left(\tau^\theta_1\otimes{1_H}\right)\circ
\tau^\theta_1\right)(\varphi_1, \dots, \varphi_p)}
\circ\left\vert\varphi_1\cdots\varphi_p\right\vert\inv[\gamma]\\
= &(-1)^p\left(\left({\bb_p}_*(\tau^\theta_1)^{\otimes p}\right)
(\varphi_1, \dots, \varphi_p)\right)[\gamma].
\endsplit
$$
See also [HS] ch.II, Theorem 3, p.126. 
This completes the proof.\qed\enddemo
\par

We conclude this section 
by describing contraction formulae deduced from the relation (2.8).
As was proved in the previous theorem for $p = 2$, 
$$
h_2 = {\bb_2}_*{h_1}^{\otimes 2} 
= \left(h_1\otimes 1_H\right)\circ \tau^\theta_1 
= \left(h_1\otimes 1_H\right)\circ h_1 \in H^2(\An; H^{\otimes 2}).
\tag 4.10
$$

We can consider other ways than $\bb_p$ of contracting the coefficients 
of the cohomology class ${h_1}^{\otimes p}$. 
For example, we may consider the $\GL(H)$-homomorphism
$$
\bb'_2: (H^*\otimes H^{\otimes 2})^{\otimes 2} \to 
H^*\otimes H^{\otimes 3}, \quad
u_{(1)}\otimes u_{(2)}\mapsto (1_H\otimes u_{(1)})\circ u_{(2)}.
$$
But we have ${\bb'_2}_*({h_1}^{\otimes 2}) = -{\bb_2}_*({h_1}^{\otimes 2}) 
= -h_2$. In fact, from (2.6), the second Johnson map $\tau^\theta_2$ gives 
the relation
$$
(h_1\otimes 1_H)\circ h_1 + (1_H\otimes h_1)\circ h_1 = 0
\in H^2(\An; H^*\otimes H^{\otimes 3}).
\tag 4.11
$$

The relation (4.11) is the IH relation among
twisted Morita-Mumford classes proposed 
by Garoufalidis and Nakamura [GN] 
and established by Morita and the author [KM2]. 
Recently Akazawa [Ak] gave an alternative proof of it
by using representation theory of the symplectic group. 
In [KM2] a more precise formula was proved. 
In \S5 we will give a simple proof of the precise  
formula using the second Johnson map $\tau^\theta_2$ 
(Theorem 5.5).\par

Let $S^0_p$ be the vertices of the Stasheff associahedron 
$K_{p+1}$ [S]. 
By definition, it is the set of all the maximal meaningful ways
of  inserting parentheses into the word $12\cdots (p+1)$ of 
$p+1$ letters. 
If we define $S^0_0 := \{1\}$, then we have 
$$
S^0_p = {\coprod}^{p-1}_{q=0} S^0_q\times S^0_{p-q-1} \tag 4.12 
$$
for $p \geq 1$. We write $\vert w\vert := p$ for $w \in S^0_p$. 
We define a sign map $\sgn: S^0_p \to \{\pm1\}$ 
by $\sgn(1) = \sgn((12)) = 1$ and 
$$
\sgn((w_1, w_2)) := (-1)^{\vert w_2\vert}\sgn(w_1)\sgn(w_2),
\quad w_1 \in S^0_q, \,\,w_2 \in S^0_{p-q-1}.
$$
Here $(w_1, w_2) \in S^0_q\times S^0_{p-q-1}$ is regarded as 
an element of $S^0_p$ by (4.12). Moreover we define a map
$$
h: S^0_p \to H^p(\An; H^*\otimes H^{\otimes (p+1)})
$$
by $h(1) := 1 = 1_H$, $h((12)) := h_1$ and
$$
h((w_1, w_2)) := (h(w_1)\otimes h(w_2))\circ h_1, \quad
w_1 \in S^0_q, \,\,w_2 \in S^0_{p-q-1}.
$$

\proclaim{Lemma 4.2} For any $w \in S^0_p$ we have 
$$
({1_H}\otimes h(w))\circ h_1 
= (-1)^{\vert w\vert}(h(w)\otimes {1_H})\circ h_1.
$$
\endproclaim
\demo{Proof} Induction on $p \geq 1$. In the case $p=1$ we have 
$({1_H}\otimes h_1)\circ h_1 = -(h_1\otimes {1_H})\circ h_1$ from (4.11). 
Suppose $p \geq 2$. For any $w = (w_1, w_2) \in S^0_q\times S^0_{p-q-1}$ 
we have, by the inductive assumption, 
$$
\split
& ({1_H}\otimes h(w))\circ h_1 = ({1_H}\otimes h(w_1)\otimes h(w_2))\circ
({1_H}\otimes h_1)\circ h_1\\
= & -({1_H}\otimes h(w_1)\otimes h(w_2))\circ (h_1\otimes {1_H})\circ h_1\\
= & -(-1)^{\vert w_1\vert}(h(w_1)\otimes {1_H}\otimes h(w_2))\circ (h_1\otimes
{1_H})\circ h_1\\
= & (-1)^{\vert w_1\vert}(h(w_1)\otimes {1_H}\otimes h(w_2))\circ ({1_H}\otimes
h_1)\circ h_1\\
= & (-1)^{\vert w_1\vert+\vert w_2\vert}(h(w_1)\otimes h(w_2)\otimes
{1_H})\circ ({1_H}\otimes h_1)\circ h_1\\
= & -(-1)^{\vert w_1\vert+\vert w_2\vert}(h(w_1)\otimes h(w_2)\otimes
{1_H})\circ (h_1\otimes {1_H})\circ h_1\\
= & (-1)^{\vert w\vert}(h(w)\otimes
{1_H})\circ h_1. \\
\endsplit
$$
This completes the induction. \qed
\enddemo

\proclaim{Proposition 4.3} We have 
$h(w) = \sgn(w) h_p \in H^p(\An; H^*\otimes H^{\otimes (p+1)})$
for any $w \in S^0_p$.
\endproclaim
\demo{Proof} We prove it by induction on $p \geq 1$. 
Immediately from the definition $h(1) = +{1_H}$ and 
$h((12)) = +h_1$. Suppose $p \geq 2$. 
For any $w = (w_1, w_2) \in S^0_q\times S^0_{p-q-1}$ 
we have, by Lemma 4.2 and the inductive assumption, 
$$
\split
& \sgn(w)h(w) = (-1)^{\vert w_2\vert}\sgn(w_1)\sgn(w_2)
(h(w_1)\otimes h(w_2))\circ h_1\\
= & (-1)^{\vert w_2\vert}\sgn(w_1)\sgn(w_2)
(h(w_1)\otimes {1_H}^{\otimes\vert w_2\vert})\circ ({1_H}\otimes
h(w_2))\circ h_1\\ = & \sgn(w_1)\sgn(w_2)
(h(w_1)\otimes {1_H}^{\otimes\vert w_2\vert})\circ
(h(w_2)\otimes {1_H})\circ h_1\\ = & 
(h_{\vert w_1\vert}\otimes {1_H}^{\otimes\vert w_2\vert})\circ
(h_{\vert w_2\vert}\otimes {1_H})\circ h_1 = h_p, \\
\endsplit
$$
which completes the induction. \qed
\enddemo

In a forthcoming paper [Ka3] we will discuss more about the   
relation between the Stasheff associahedron $K_{p+1}$ and 
the cohomology class $h_p$.

\beginsection 5. Mapping Class Groups.\par

Let $g \geq 1$ be a positive integer, 
$\Sigma_{g, 1}$ a $2$-dimensional oriented compact connected 
$C^\infty$ manifold of genus $g$ with $1$ boundary component. 
We choose a basepoint $*$ on the boundary $\partial \Sigma_{g, 1}$. 
The fundamental group $\pi_1(\Sigma_{g, 1}, *)$ is a free group of 
rank $2g$. 
Taking a symplectic generator system $\{x_i\}^{2g}_{i=1} \subset
\pi_1(\Sigma_{g, 1}, *)$, we identify $\pi_1(\Sigma_{g, 1}, *) = 
F_{2g}$. This induces a natural isomorphism
$$
H = H_1(F_{2g}; R) = H_1(\Sigma_{g, 1}; R).
$$
A simple loop parallel to the boundary gives a word
$$
w_0 := {\prod}^g_{i=1} x_ix_{g+i}{x_i}\inv{x_{g+i}}\inv.
$$
The intersection number $\cdot$ on the surface $\Sigma_{g, 1}$ satisfies
$$
X_i\cdot X_{g+j} = \delta_{i, j} \quad\text{and}\quad
X_i\cdot X_j = X_{g+i}\cdot X_{g+j} = 0
$$
for $1 \leq i, j \leq g$, and the intersection form $I$ is given by 
$$
I = {\sum}^g_{i=1}(X_i\otimes X_{g+i} -  X_{g+i}\otimes X_i)
\in H^{\otimes 2}. 
$$
We denote by  $\{\xi_i\}^{2g}_{i=1} \subset H^*$ the dual basis of 
$\{X_i\}^{2g}_{i=1} \subset H$. 
The Poincar\'e duality $\vartheta :=\cap [\Sigma_{g, 1}]:  H^* = H^1(\Sigma_{g,
1}, \partial \Sigma_{g, 1}; R) 
\to H_1(\Sigma_{g, 1};\allowmathbreak R) = H$ is, by definition, 
the cap product by the fundamental class $[\Sigma_{g, 1}] \in H_2(\Sigma_{g,
1}, \partial \Sigma_{g, 1}; \bZ)$. Then we have 
$$
\split
1 & = X_i\cdot X_{g+i} 
= \langle \vartheta\inv(X_i)\cup\vartheta\inv(X_{g+i}), 
[\Sigma_{g, 1}]\rangle\\
& = \langle \vartheta\inv(X_i), \vartheta\inv(X_{g+i})\cap 
[\Sigma_{g, 1}]\rangle 
= \langle \vartheta\inv(X_i), X_{g+i}\rangle,
\endsplit
$$
and so on. Hence we obtain
$$
\vartheta\inv(X_i) = \xi_{g+i} = X_i\cdot\empty, \quad\text{and}\quad
\vartheta\inv(X_{g+i}) = -\xi_{i} = X_{g+i}\cdot\empty
$$
for $1 \leq i \leq g$, or equivalently, 
the Poincar\'e duality $\vartheta$ and its inverse $\vartheta\inv$ are 
given by 
$$
\aligned
& \vartheta: H^* \overset\cong\to\to H, \quad
\ell \mapsto -(\ell\otimes 1_H)(I) = (1_H\otimes\ell)(I), \quad \text{and}\\
& \vartheta\inv: H \overset\cong\to\to H^*, \quad
Y \mapsto Y\cdot\empty.
\endaligned\tag 5.1
$$ 
In this section we identify 
$H^* = H$ by the Poincar\'e duality (5.1), which is equivariant under the
action of the mapping class group 
$\M1g := \pi_0\operatorname{Diff}(\Sigma_{g, 1} 
\operatorname{rel} \partial \Sigma_{g, 1})$. \par
There is another sign convention on cap products, 
which our previous papers [KM1] [KM2] followed. 
In that convention we have 
$$
\langle \vartheta\inv(X_i)\cup\vartheta\inv(X_{g+i}), 
[\Sigma_{g, 1}]\rangle 
= \langle \vartheta\inv(X_{g+i}), \vartheta\inv(X_i)\cap 
[\Sigma_{g, 1}]\rangle,
$$
and so $\vartheta(\ell) = (\ell\otimes 1_H)(I)$ and $\vartheta\inv(Y) = 
-Y\cdot\empty$. Then we obtain 
$(p+2)!h_p\vert_{\M1g} = -m_{0, p+2}$ and $p!\h_p\vert_{\M1g} = m_{1, p}$ in
(5.8) and (5.18). \par
The mapping class group $\M1g$ acts on the fundamental group 
$\pi_1(\Sigma_{g, 1}, *) = F_{2g}$. This induces a group homomorphism 
$\M1g \to \Aut(F_{2g})$. We prove the pull-back of the cohomology classes 
$h_p$ and $\h_p$ to the group $\M1g$ are twisted Morita-Mumford classes 
[Ka1] (Theorem 5.1 and Corollary 5.4). 
Furthermore we give a simple proof of a precise version [KM2]
of the IH-relation [GN]. \par
We introduce some variants of the mapping class group $\M1g$ 
in order to recall the definition of the Morita-Mumford classes 
[Mo1] [Mu] and that of the twisted ones [Ka1] [KM1]. 
Collapsing the boundary $\partial \Sigma_{g, 1}$ into a single point $*$, 
we obtain a $2$-dimensional oriented closed connected $C^\infty$ 
manifold of genus $g$, $\Sigma_g = \Sigma_{g, 1}/\partial \Sigma_{g, 1}$. 
We write simply $\pi_1 := \pi_1(\Sigma_{g}, *)$ and $\pi^0_1 := \pi_1
(\Sigma_{g, 1}, *)$. Let $N_g$ be the kernel of the collapsing 
homomorphism $c: \pi^0_1 \to \pi_1$. Then we have a group extension 
$$
N_g \to \pi^0_1 \overset{c}\to\to \pi_1.\tag 5.2
$$
Let $\Mg$ and $\Mgstar$ be the mapping class groups for the 
surface $\Sigma_g$ and the pointed one $(\Sigma_g, *)$, 
respectively, that is, $\Mg := \pi_0\operatorname{Diff}_+(\Sigma_g)$ and 
$\Mgstar := \pi_0\operatorname{Diff}_+(\Sigma_g, *)$. Forgetting the 
basepoint induces a group extension
$$
\pi_1 \to \Mgstar \overset\pi\to\to \Mg, \tag 5.3
$$
and collapsing the boundary a central extension of groups
$$
\bZ \to \M1g \overset\varpi\to\to \Mgstar. \tag 5.4
$$
The kernel $\Ker\,\varpi \cong \bZ$ is generated by the Dehn twist 
along a simple loop parallel to the boundary, 
which acts on the fundamental group $\pi^0_1$ 
by conjugation by the word $w_0$.\par

We denote the Euler class of the central extension (5.4) by 
$e:= \operatorname{Euler}(\bZ \to \M1g \overset\varpi\to\to \Mgstar) 
\in H^2(\Mgstar; \bZ)$. The extension (5.3) gives the Gysin map 
$$
\pi_!: H^*(\Mgstar; M) \to H^{*-2}(\Mg; M)
$$
for any $\Mg$-module $M$. 
The Morita-Mumford class $e_i$, $i \geq 1$, is defined to be 
the Gysin image of the $(i+1)$-st power of the Euler class $e$
$$
e_i := \pi_!(e^{i+1}) \in H^{2i}(\Mg; \bZ).
$$
The fiber product $\M1g\times_\Mg\Mgstar$ is identified with the 
semi-direct product $\pi_1\rtimes\M1g$ by the isomorphism  
$(\varphi, \psi) \in \M1g\times_\Mg\Mgstar\mapsto (\psi\varpi(\varphi)\inv, 
\varphi) \in \pi_1\rtimes\M1g$. 
We denote the first and the second projections of the product 
 $\M1g\times_\Mg\Mgstar$ by $\pi: \pi_1\rtimes\M1g \to \M1g$ and 
$\overline{\pi}: \pi_1\rtimes\M1g \to \Mgstar$, respectively, 
and write $\ebar := \overline{\pi}^*(e) \in H^2(\pi_1\rtimes\M1g ;\bZ)$. 
Then the pullback $e_i = \varpi^*\pi^*e_i \in H^{2i}(\M1g; \bZ)$ is given by 
$$
e_i = \pi_!(\ebar^{i+1}) \in H^{2i}(\M1g;\bZ). \tag 5.5
$$
As in \S4, we can consider the $1$-cocycle $k_0: \pi_1\rtimes\M1g \to H$, 
$(\gamma, \varphi) \mapsto [\gamma]$. 
The twisted Morita-Mumford class $m_{i, j}$, $i, j \geq 0$, $2i+j \geq 2$, 
is defined to be the Gysin image of $\ebar^i{k_0}^j$
$$
m_{i, j} := \pi_!(\ebar^i{k_0}^j) \in H^{2i+j-2}(\M1g; \Lambda^jH).
\tag 5.6
$$
Here ${k_0}^j$ is the $j$-th exterior power of $k_0$, and so we have 
$$
{k_0}^j = j!{k_0}^{\otimes j}\in H^p(\M1g; H^{\otimes j}).
\tag 5.7
$$
Now we have

\proclaim{Theorem 5.1}
$$
h_p = \pi_!\left({k_0}^{\otimes (p+2)}\right) 
\in H^p(\M1g; H^{\otimes (p+2)}) 
$$
for each $p \geq 1$.
\endproclaim

It is an immediate consequence of Lemma 5.4 in [KM2]. 
But we will give a self-contained proof of the theorem. 
From (5.6) and (5.7) we obtain
$$
(p+2)!h_p = m_{0, p+2}.\tag 5.8
$$
\par
We should remark on the identification of $\Hom(H, H^{\otimes (p+1)})$ 
with $H^{\otimes (p+2)}$. We denote it by 
$$
t_p: \Hom(H, H^{\otimes (p+1)}) = H^*\otimes H^{\otimes (p+1)}
\overset{\vartheta\otimes 1}\to\longrightarrow H^{\otimes
(p+2)}.
$$
From (5.1) we have 
$$
t_p(u) = -(1,2,\dots, p+1, p+2)(u\otimes 1_H)(I)
\tag 5.9
$$
for any $u \in \Hom(H, H^{\otimes (p+1)})$. 
Here the cyclic permutation $(1,2,\dots, p+1, p+2)$ acts on 
$H^{\otimes(p+2)}$ by permuting the components of the tensors 
in $H^{\otimes(p+2)}$.
\par
In order to prove the theorem we construct a cohomology class 
introduced in [Mo2]
$$
\nu \in H^2(\pi_1\rtimes\M1g; \bZ)
$$
in an algebraic way similar to [KM2]. 
In this section we write simply
$$
\Mbar := \pi_1\rtimes \M1g \quad\text{and}\quad \Mzbar := \pi^0_1\rtimes \M1g.
$$
The collapsing homomorphism gives a group extension
$$
N_g \to \Mzbar \overset{c}\to\to \Mbar. \tag 5.10
$$
The Lyndon-Hochschild-Serre spectral sequence of the extension (5.2) 
gives 
$$
H^p(\pi_1; H^1(N_g; \bZ)) = 0, \quad\text{if $p \geq 1$,}\tag 5.11
$$
and an $\M1g$-invariant isomorphism
$$
d_2: H^1(N_g; \bZ)^{\pi_1}\overset\cong\to\to 
H^2(\Sigma_g; \bZ) \cong \bZ. 
$$
Choose a Magnus expansion $\theta \in \Theta_{2g, \bZ}$. 
We have $\theta_2(w_0) = I$. 
From (3.2) for $m=2$ follows
$$
\theta_2(\gamma w_0\gamma\inv) 
= \theta_2(\gamma w_0\gamma\inv{w_0}\inv w_0)
= \theta_2(\gamma w_0\gamma\inv{w_0}\inv) + \theta_2(w_0)
= \theta_2(w_0)
$$
for any $\gamma \in \pi^0_1$.
Since $N_g$ is the normal closure of 
the word $w_0$, we can define an $\Mbar$-invariant homomorphism 
$\nu_0: N_g \to \bZ$ by 
$$
\nu_0(\delta)I = -\theta_2(\delta) \in H^{\otimes 2}
$$
for $\delta \in N_g$.
Consider the transgression of the Lyndon-Hochschild-Serre spectral 
sequence of the extension (5.10)
$$
d_2: H^0(\Mbar; H^1(N_g; \bZ)) \to H^2(\Mbar; \bZ).
$$
If we choose a map $\widehat{\empty}: \pi_1 \to \pi^0_1$, 
$\gamma \mapsto \widehat{\gamma}$, satisfying 
$c(\widehat{\gamma}) = \gamma$ for any $\gamma \in \pi_1$ and 
$\widehat{1} = 1$, 
then the $2$-cochain $\widehat{\nu}$ defined by 
$$
\widehat{\nu}(\gamma_1\varphi_1, \gamma_2\varphi_2) := 
\nu_0(\widehat{\gamma_1\varphi_1(\gamma_2)}
\varphi_1(\widehat{\gamma_2})\inv \widehat{\gamma_1}\inv)
\tag 5.12
$$
for $\gamma_1\varphi_1$ and $\gamma_2\varphi_2 \in \Mbar$ 
represents $$
\nu := d_2\nu_0 \in H^2(\Mbar; \bZ).
$$ 
Moreover we have
$$
\pi_!(\nu) = \langle\widehat{\nu}, [\Sigma_g]\rangle = 1,
\tag 5.13
$$
where $[\Sigma_g] \in H_2(\Sigma_g; \bZ)$ is the fundamental class. \par
To prove (5.13) we consider the homomorphism $\phi: \pi^0_1 \to F_2$ 
given by $\phi(x_1) = x_1$, $\phi(x_{g+1}) = x_2$ and  $\phi(x_j) =
\phi(x_{g+j}) = 1$ for $j \geq 2$, which induces a homomorphism 
of group extensions
$$
\CD
N_g  @>>> \pi^0_1 @>{c}>> \pi_1\\
@V{\phi}VV @V{\phi}VV @V{\phi}VV \\  
[F_2, F_2] @>>> {F_2} @>>> {F_2}\abel.
\endCD\tag 5.14
$$
We regard ${F_2}\abel$ as the fundamental group of 
the $2$-dimensional (real) torus $T^2$. 
Then $\phi$ preserves the orientations, that is, 
$\phi_*[\Sigma_g] = [T^2] \in H_2(T^2; \bZ)$. 
The fundamental class is given by a normalized bar $2$-chain
$$
[x_1\vert x_2] + [x_1x_2\vert {x_1}\inv] - [{x_1}\vert{x_1}\inv].
$$
See, e.g., [Me] p\.245. If we use the map $\widehat{\empty}: 
{F_2}\abel \to F_2$, ${x_1}^a{x_2}^b \mapsto {x_1}^a{x_2}^b$, 
then we have $\langle \widehat{\nu}, [T^2]\rangle 
= \nu_0(x_2x_1{x_2}\inv{x_1}\inv) = 1$. 
From (5.14) we have $\langle\widehat{\nu}, [\Sigma_g]\rangle = 
\langle\widehat{\nu}, \phi_*[\Sigma_g]\rangle = 
\langle\widehat{\nu}, [T^2]\rangle = 1$ for any $g \geq 1$. 
This proves (5.13).\qed\par
\medskip
The following is the key to the proof of Theorem 5.1.

\proclaim{Lemma 5.2} $\nu k_0 = 0 \in H^3(\Mbar; H)$.
\endproclaim
\demo{Proof} The lemma is an immediate consequence of Theorem 5.1 (ii) 
in [KM2]. But we give a self-contained proof of it. 
The Lyndon-Hochschild-Serre spectral sequence of the semi-direct product 
$\Mbar = \pi_1\rtimes\M1g$ gives an isomorphism 
$$
\pi^*: H^1(\M1g; H^1(N_g)^{\pi_1}\otimes H) \overset\cong\to\to 
H^1(\Mbar; H^1(N_g)\otimes H), \tag 5.15
$$
since $H^1(\pi_1; H^1(N_g)) = 0$ (5.11). 
Consider the homomorphism $s: \M1g \to \Mbar$, 
$\varphi \mapsto (1, \varphi)$. Then
$$
s^*: H^1(\Mbar; H^1(N_g)\otimes H) \to H^1(\M1g; H^1(N_g)\otimes H)
$$
is an isomorphism. In fact, $\pi^*: H^1(\M1g; H^1(N_g)\otimes H) 
\to H^1(\Mbar; H^1(N_g)\otimes H)$ is surjective from the isomorphism 
(5.15). Clearly we have $s^*\pi^* = 1$ on 
$H^1(\M1g; H^1(N_g)\otimes H)$. Hence $s^*$ is the inverse of the 
isomorphism $\pi^*$.\par
Since $s^*k_0 = 0$, we have $\nu_0k_0 = \pi^*s^*(\nu_0k_0) = 0 
\in H^1(\Mbar; H^1(N_g)\otimes H)$. 
Consequenctly $\nu k_0 = d_2(\nu_0k_0) = 0$, as was to be shown.\qed
\enddemo

\demo{Proof of Theorem 5.1}
We denote the Gysin map of the semi-direct product
$\Mzbar = \pi^0_1\rtimes\M1g$ by $\pi_\sharp: H^*(\Mzbar; M) 
\to H^{*-1}(\M1g; H\otimes M)$ for any $\M1g$-module $M$. 
From the definition of $h_p$, we have 
$h_p = \pi_\sharp({k_0}^{\otimes (p+1)})$. 
Consider the $1$-cochain $\widehat{\theta_2} \in 
C^1(\Mbar; H^{\otimes 2})$ defined by 
$$
\widehat{\theta_2}(\gamma\varphi) := \theta_2(\widehat{\gamma})
$$
for $\gamma\varphi \in \pi_1\rtimes\M1g = \Mbar$. Then we have
$$
\split
& (d\widehat{\theta_2})(\gamma_1\varphi_1, \gamma_2\varphi_2) 
= \varphi_1\theta_2(\widehat{\gamma_2}) 
- \theta_2(\widehat{\gamma_1\varphi_1(\gamma_2)}) 
+ \theta_2(\widehat{\gamma_1})\\
=\,& \varphi_1\theta_2(\widehat{\gamma_2}) 
- \theta_2(\widehat{\gamma_1\varphi_1(\gamma_2)}
\varphi_1(\widehat{\gamma_2})\inv\widehat{\gamma_1}\inv
\widehat{\gamma_1}\varphi_1(\widehat{\gamma_2}))
+ \theta_2(\widehat{\gamma_1})\\
=\,& \varphi_1\theta_2(\widehat{\gamma_2}) 
- \theta_2(\varphi_1(\widehat{\gamma_2})) 
+ \theta_2(\varphi_1(\widehat{\gamma_2})) 
- \theta_2(\widehat{\gamma_1\varphi_1(\gamma_2)}
\varphi_1(\widehat{\gamma_2})\inv\widehat{\gamma_1}\inv) \\
& \qquad\qquad\qquad\qquad\qquad\qquad\qquad
-\theta_2(\widehat{\gamma_1}\varphi_1(\widehat{\gamma_2})) +
\theta_2(\widehat{\gamma_1})\\
=\,& -\tau^\theta_1(\varphi_1)\vert\varphi_1\vert[\gamma_2]
- [\gamma_1]\otimes \vert\varphi_1\vert[\gamma_2]
+ \widehat{\nu}(\gamma_1\varphi_1, \gamma_2\varphi_2)I
\endsplit
$$
for any $\gamma_1\varphi_1$ and $\gamma_2\varphi_2 \in \Mbar$.
This means
$$
\tau^\theta_1\circ k_0 + {k_0}^{\otimes 2} = \nu I 
\in H^2(\Mbar; H^{\otimes 2}).\tag 5.16
$$
Hence we have 
$$
\pi_!({k_0}^{\otimes 2}) = I \in H^0(\M1g; H^{\otimes 2})
\tag 5.17
$$
from (5.13). 
Moreover, from Lemma 5.2 and (5.16), we have
$$
(h_1\otimes {1_H}^{\otimes p})\circ {k_0}^{\otimes (p+1)}
+ {k_0}^{\otimes (p+2)} = 0 
$$
for each $p \geq 1$.
If we denote $h'_p := \pi_!({k_0}^{\otimes (p+2)}) \in H^p(\M1g; H^{\otimes
(p+2)})$, then $(h_1\otimes {1_H}^{\otimes p})\circ h'_{p-1} + h'_p = 0$, 
and so
$$
\align
h'_p = & (-1)^p (h_1\otimes{1_H}^{\otimes p})\circ(h_1\otimes{1_H}^{\otimes
(p-1)})\circ\cdots\circ (h_1\otimes{1_H})\circ h'_0,\\
= & (-1)^p(h_p\otimes 1_H)(I)
\endalign
$$
from Theorem 4.1 and (5.17). 
Consequently, from (5.9) and the commutativity of the cup product, 
we obtain
$$
\split
{t_p}_*(h_p) = & -(1,2,\dots, p+1, p+2)_*(h_p\otimes 1_H)(I)\\
= & (-1)^{p+1}(1,2,\dots, p+1, p+2)_*h'_p\\
= & (-1)^{p+1}(1,2,\dots, p+1, p+2)_*\pi_!({k_0}^{\otimes (p+2)})\\
= & \pi_!({k_0}^{\otimes (p+2)}), 
\endsplit
$$
as was to be shown.
\qed 
\enddemo

Next we study the cohomology class $\h_p$. 
We denote by $\mu: H^{\otimes 2}\to \bZ$ the intersection product 
on the surface $\Sigma_{g, 1}$. 
Recall the following theorem due to Morita.

\proclaim{Theorem 5.3}{\rm (Morita [Mo2], Theorem 1.3.)} $$
\mu_*({k_0}^{\otimes
2}) = 2\nu - \ebar 
\in H^2(\pi_1\rtimes\M1g; \bZ).
$$
\endproclaim 

An algebraic proof of it is given in [KM2], Theorem 6.1.\par
From Lemma 5.2 we have $\mu_*({k_0}^{\otimes 2})\otimes 
{k_0}^{\otimes p} = -\ebar\otimes{k_0}^{\otimes p}$. 
Theorem 5.1 implies
$$
\h_p = (\mu\otimes {1_H}^{\otimes p})_*h_p = (\mu\otimes {1_H}^{\otimes
p})_*\pi_!({k_0}^{\otimes (p+2)}) = -\pi_!(\ebar\otimes{k_0}^{\otimes p}).
$$
Hence we obtain
\proclaim{Corollary 5.4} $$
\h_p = -\pi_!(\ebar\otimes{k_0}^{\otimes p})
\in H^p(\M1g; H^{\otimes p}).
$$  
\endproclaim
From (5.7) follows
$$
p!\h_p = -m_{1, p}. \tag 5.18
$$
If we contract the coefficients of $h_p$ by an iteration of the product $\mu$, 
then we obtain the (original) Morita-Mumford class $e_i$. See [KM1] [KM2].
Hence each of the $e_i$'s is given by a certain algebraic combination of 
copies of $h_1$. See also Theorems 4.1 and 5.3. So we may consider the
cohomology class $h_1$ as ``the unique elementary particle" for all the 
Morita-Mumford classes.\par

Finally we study a consequence of the relation (2.6) on the mapping class 
group $\Mgstar$. 
The conjugation by the word $w_0$, $\iota(w_0)$, corresponds to 
a generator of the kernel $\Ker(\varpi: \M1g\to \Mgstar)$ and we have 
$\tau^\theta_1\iota(w_0) = 0$ from (2.11). 
Hence the first Johnson map $\tau^\theta_1$ can be regarded as a $1$-cocycle 
$\tau^\theta_1$ on the group $\Mgstar$. 
We denote $h_1 := [\tau^\theta_1] \in H^1(\Mgstar; H^*\otimes H^{\otimes 2})$. 
The following is a precise version of the IH-relation [GN].
\proclaim{Theorem 5.5}{\rm ([KM2], Theorem 1.3 (iii))}
$$
(h_1\otimes 1_H)\circ h_1 + (1_H\otimes h_1) \circ h_1 
= e(I\otimes 1_H - 1_H\otimes I) \in H^2(\Mgstar; H^*\otimes H^{\otimes 3}).
$$
\endproclaim
\demo{Proof} From (2.6) we have 
$$
-d\tau^\theta_2 = (\tau^\theta_1\otimes 1_H + 1_H\otimes\tau^\theta_1) \cup 
\tau^\theta_1
$$
on the group $\M1g$. 
The Gysin sequence of the extension (5.4) 
$$
H^0(\Mgstar; H^*\otimes H^{\otimes 3}) \overset{\cup e}\to\to
H^2(\Mgstar; H^*\otimes H^{\otimes 3}) \overset{\varpi^*}\to\to
H^2(\M1g; H^*\otimes H^{\otimes 3}) 
$$
implies 
$$
(h_1\otimes 1_H)\circ h_1 + (1_H\otimes h_1) \circ h_1 
= e \tau^\theta_2(\iota(w_0)) = e(I\otimes 1_H - 1_H\otimes I)
$$
from (2.12). This proves the theorem.\qed
\enddemo

As was stated above, the coefficients $H^*\otimes H^{\otimes 3}$ in the 
theorem are identified with $H^{\otimes 4}$ by the map $t_3$ in (5.9). 
In [KM2] a closed trivalent graph describes an $Sp_{2g}(\bQ)$-invariant 
of the algebra $H^*(\La^3H_\bZ; \bQ) \cong \La^*(\La^3H_\bQ)$. 
In this context the cohomology class $h_1 = [\tau^\theta_1] \in 
H^1(\M1g; H^{\otimes 3})$ corresponds to the open star of each vertex 
on the graph. 
A subgraph shaped like the letter H means the twisted cohomology class 
$(1_H\otimes h_1)\circ h_1$, while one like the letter I means 
$-(h_1\otimes 1_H)\circ h_1$ because of the commutativity of 
the cup product. 
Similarly $I\otimes 1_H$ and $1_H\otimes I$ are interpreted as 
suitable edges in $\Gamma_1\setminus \tau_1$ and $\Gamma_2\setminus \tau_2$ 
in [KM2], respectively. 
Hence Theorem 1.3 (iii) in [KM2] follows from our theorem. 

\beginsection 6. The Abelianization of $\IAn$.\par

The group $\IAn$ is defined to be the kernel of the homomorphism
$\vert\,\,\vert: \An \to \GL(\HZ)$ induced by the abelianization 
$\Fn \to \Fn\abel = \HZ$. 
In other words, $\IAn = A(1)$ in \S3.
Classically it is called {\it the induced automorphism  group}. 
In this section we compute the abelianization of $\IAn$  
by evaluating the first Johnson map on the generators 
of the group $\IAn$ given by Magnus [M2],  
and give some consequences of the computation. 
Here it should be  remarked that
S\. Andreadakis has already studied the abelianization 
$\IAn\abel$ in [An]. \par

First recall the second exterior power $\Lambda^2\HZ$ of 
$\HZ = \Fn\abel$. Let $S^2(\HZ) \subset \HZ^{\otimes 2}$ be 
the $\bZ$-submodule generated by the set $\{Y\otimes Y; \,\,
Y \in \HZ\}$. By definition we have $\Lambda^2\HZ = 
\HZ^{\otimes 2}/S^2(\HZ)$. We define an homomorphism 
$\alpha_2: \HZ^{\otimes 2}\to \HZ^{\otimes 2}$ by 
$\alpha_2(X_i\otimes X_j) = X_i\otimes X_j - X_j\otimes X_i$ 
for $1\leq i, j\leq n$. We have $\Ker\,\alpha_2 = S^2(\HZ)$.
This induces an injective homomorphism $\overline{\alpha_2}: 
\Lambda^2\HZ \to \HZ^{\otimes 2}$. Throughout this section 
we regard $\Lambda^2\HZ$ as a submodule of $\HZ^{\otimes 2}$ 
by the injection $\overline{\alpha_2}$.\par

Now fix a $\bZ$-valued Magnus expansion $\theta\in \Theta_{n, \bZ}$. 
It gives the first Johnson map $\tau^\theta_1: \An \to {{H_{\Bbb
Z}}^*\otimes{H_{\Bbb Z}}^{\otimes 2}}$. 
As was proved in Theorem 3.1, 
the restriction of $\tau^\theta_1$ to $\IA_n = A(1)$ is 
equal to the first Johnson homomorphism $\tau_1$, 
which is independent of the choice of $\theta$. 
\proclaim{Theorem 6.1}
The first Johnson homomorphism $\tau_1$ induces an isomorphism
$$
\tau_1: \IAn\abel \overset{\cong}\to\to \tensors,
$$
which is equivariant under the action of
$\GL(\HZ) = \An/\IAn$. Especially the abelianization 
$\IAn\abel$ is free abelian of rank $n^2(n-1)/2$, and 
the commutator subgroup of $\IAn$ coincides 
with $\Ker\,\tau_1$
$$
[\IAn, \IAn] = \Ker\,\tau_1 = A(2). \tag 6.1
$$
\endproclaim

Andreadakis [An] proved the theorem for the case $n=3$. 
All we need to prove it are due to W\. Magnus. 
So it had been likely proved by someone contemporary with 
Magnus or Andreadakis. Comparing it with 
Johnson's result computing the abelianization of 
the Torelli groups [J3], the reader would find how simpler 
the automorphism groups of free groups are than the 
mapping class groups for surfaces.

\demo{Proof of Theorem 6.1} 
According to Magnus [M2], 
the group $IA_n$ is generated by the following 
automorphisms
$$
\align
& K_{i, l}: x_i \mapsto x_lx_i{x_l}^{-1}, 
\quad x_j \mapsto x_j \quad (j \neq i)\\
& K_{i, l, s}: x_i \mapsto x_ix_lx_s{x_l}^{-1}{x_s}^{-1}, 
\quad x_j \mapsto x_j \quad (j \neq i).
\endalign
$$
Here the indices run over the sets 
$\{(i, l); 1 \leq i, l \leq n, \,\, i \neq l\}$ and 
$\{(i, l, s); 1 \leq i, l, s \leq n, \,\, i 
\neq l < s \neq i\}$ respectively. 
The number of the generators is 
$n(n-1) + n\frac{1}{2}(n-1)(n-2) = \frac{1}{2}n^2(n-1)$.
Hence we have a surjection $p_n: \bZ^{n^2(n-1)/2} \to 
\IAn\abel$. \par
Now we denote by $\theta_2: \Fn \to \HZ^{\otimes 2}$ 
the second component of the expansion $\theta$ 
as before. From (3.2) for $m=1$ we have 
$$
\theta_2(\gamma\delta\gamma\inv\delta\inv) =
[\gamma][\delta] - [\delta][\gamma] 
\tag 6.2
$$
for any $\gamma$ and $\delta \in \Fn$.
Since $\theta_2(\gamma\delta\gamma\inv) =
\theta_2(\gamma\delta\gamma\inv\delta\inv\delta) = 
\theta_2(\gamma\delta\gamma\inv\delta\inv) + \theta_2(\delta)$, 
we have 
$$
\theta_2(\gamma\delta\gamma\inv) - \theta_2(\delta) 
= [\gamma][\delta] - [\delta][\gamma].
\tag 6.3
$$
\par

Let $\{\ell_i\}^n_{i=1} \subset \HZ^*$ be the dual basis of 
the basis $\{X_i\}^n_{i=1} \subset \HZ$. From (2.7), for any
$\varphi 
\in \IAn$ and $\gamma\in \Fn$, we have 
$
\tau^\theta_1(\varphi)[\gamma] = 
\theta_2(\varphi(\gamma)) - \theta_2(\gamma).
$
From (6.2) and (6.3) we obtain
$$
\aligned
\tau^\theta_1(K_{i, l}) &= 
\ell_i\otimes\tau^\theta_1(K_{i,l})(X_i) = 
\ell_i\otimes\left(\theta_2(x_lx_i{x_l}\inv) -
\theta_2(x_i)\right)\\ & = 
\ell_i\otimes (X_l X_i - X_i X_l)\\
\tau^\theta_1(K_{i, l, s}) &= 
\ell_i\otimes\tau^\theta_1(K_{i,l,s})(X_i) = 
\ell_i\otimes\left(\theta_2(x_ix_lx_s{x_l}\inv{x_s}\inv) -
\theta_2(x_i)\right)\\ & = 
\ell_i\otimes (X_l X_s - X_s X_l).
\endaligned\tag 6.4
$$
These form exactly a $\bZ$-free basis of $\tensors$. 
Therefore the composite
$$
\tau^\theta_1\circ p_n: \bZ^{n^2(n-1)/2}
\overset{p_n}\to\longrightarrow
\IAn\abel \overset{\tau^\theta_1}\to\longrightarrow \tensors
$$
is an isomorphism. Since $p_n$ is surjective, 
the homomorphism $\tau_1 = \tau^\theta_1: \IAn\abel \to \tensors$ 
is an isomorphism. \par
The isomorphism $\tau^\theta_1: \IAn\abel \to \tensors$
is $\GL(\HZ)$-equivariant, because we have
$$
\split
 \tau^\theta_1(\varphi\psi\varphi\inv)
= &\,\tau^\theta_1(\varphi\psi\varphi\inv) 
- \tau^\theta_1(\varphi\varphi\inv)\\
= &\, \vert\varphi\vert\tau^\theta_1(\psi) 
+ \vert\varphi\psi\vert\tau^\theta_1(\varphi\inv)
- \vert\varphi\vert\tau^\theta_1(\varphi\inv)
= \vert\varphi\vert\tau^\theta_1(\psi) 
\endsplit
$$
for any $\varphi \in \An$ and $\psi \in \IAn$.
This completes the proof of Theorem 6.1.
\qed\enddemo

We define the group $\IOn$ to be the kernel of the homomorphism 
$\vert\,\,\vert: \On \to \GL(\HZ)$ induced
by the abelianization. We have $\IOn = \IAn/\Inn$. 
Here $\Inn$ is, by definition, the image of 
the homomorphism $\iota$ in (2.9).
From (2.11) the induced homomorphism
$$
\iota_*: \Fn\abel = \HZ \to 
\IAn\abel \overset{\tau_1}\to\cong \tensors
$$
is given by 
$$
\iota_*(Y)Z = YZ - ZY \in \La^2\HZ \tag 6.5
$$
for $Y, Z \in \HZ$, which is an injection, and whose image 
is a direct summand of $\tensors$ as a $\bZ$-module. 
Hence
\proclaim{Theorem 6.2} We have a $\GL(\HZ)$-equivariant 
isomorphism 
$$
\IOn\abel \cong (\tensors)/\iota_*(\HZ),
$$
where $\iota_*$ is the homomorphism given in (6.5).
Especially the abelianization $\IOn\abel$ is free abelian 
of rank $(n+1)n(n-2)/2$.
\endproclaim

Let $q$ be a prime integer. 
The congruence IA-automorphism group $IA_{n, q}$ is 
defined to be the kernel of the natural homomorphism 
$\An \to \GL(H_{\bZ/q})$. 
Using the first Johnson map $\tau^\theta_1$ 
for a $\bZ/q$-valued Magnus expansion $\theta$, 
T\. Satoh [Sa2] computes the abelianization of $IA_{n, q}$.
\par

Theorem 6.1 has an application to twisted cohomology 
of the group $\An$ with values in a $GL(\HZ)$-module. 
Let $\overline{\Gamma}$ be a subgroup of $\GL(\HZ)$, 
$\Gamma \subset \An$ the preimage of $\overline{\Gamma}$,
and $M$ a $\bZ[\frac12][\overline{\Gamma}]$-module. 
We denote by $\pi: \Gamma  \to \overline{\Gamma}$ 
the natural projection, and by $\jmath$ the inclusion 
$\jmath: \IAn \hookrightarrow \Gamma$. 
We have the Lyndon-Hochschild-Serre  spectral sequence
$$
E^{p, q}_2 = E^{p, q}_2(M) = 
H^p(\overline{\Gamma}; H^q(\IAn; M)) 
\Rightarrow H^{p+q}(\Gamma; M) \tag 6.6
$$
of the group extension $\IAn \overset\jmath\to\to \Gamma
\overset\pi\to\to
\overline{\Gamma}$.

\proclaim{Proposition 6.3} For any $p \geq 0$ we have
$$
d^{p,1}_2 = 0: H^p(\overline{\Gamma}; H^1(\IAn; M)) 
\to H^{p+2}(\overline{\Gamma}; M).
$$
In the case $p = 0$ we have a natural decomposition 
$$
H^1(\Gamma; M) = \Hom(\tensors, M)^{\overline{\Gamma}}\oplus
H^1(\overline{\Gamma}; M).
$$
\endproclaim
\demo{Proof} We denote by $\varpi_2: \HZ^{\otimes 2} \to \La^2
\HZ$ the natural projection. Twice the abelianization 
$\tau^\theta_1: \IAn\to \tensors$ extends to the crossed homomorphism 
$$
\widetilde{\tau}:= (1\otimes \varpi_2)\circ \tau^\theta_1:
\An \to \HZ^*\otimes \HZ^{\otimes 2} \to \tensors 
$$
defined on the whole automorphism group $\An$. 
Since the action of the subgroup $\IAn$ on the module $M$ is 
trivial, we have the Kronecker product
$$
\kappa: \tensors\otimes H^1(\IAn; M) = H_1(\IAn)\otimes H^1(\IAn;
M) \to M.
$$ 
For any $p$-cocycle $f \in Z^p(\overline{\Gamma}; H^1(\IAn; M))$ 
the $(p+1)$-cocycle
$$
\kappa_*(\widetilde{\tau}\cup\pi^*(\frac12 f)) 
\in Z^{p+1}({\Gamma}; M)
$$
is an element of the $p$-th filter $A^{p+1}\cap A^*_p$ in [HS] ch\.II,
p\.119.  It is clear that the cocycle 
$\kappa_*(\widetilde{\tau}\cup\pi^*(\frac12 f))$
induces the cocyle $f$ in $E^{*\,p,1}_1 = C^p(\overline{\Gamma}; 
H^1(\IAn; M))$. Thus the cocycle $f$ extends to a cocycle 
defined on the whole $\Gamma$, so that $d^{p, 1}_2[f] = 0 
\in E^{p+2, 0}_2$.
In the case $p = 0$ we have an exact sequence
$$
0 \to H^1(\overline{\Gamma}; M) \overset{\pi^*}\to\to
H^1(\Gamma; M) \overset{\jmath^*}\to\to
H^0(\overline{\Gamma}; H^1(\IAn; M)).
$$
The homomorphism
$$
[f] \in H^0(\overline{\Gamma}; H^1(\IAn; M)) \mapsto
[\kappa_*(\widetilde{\tau}\cup\pi^*(\frac12 f))] \in 
H^1(\Gamma; M)
$$
is a right inverse of the homomorphism $\jmath^*$. Hence we have 
$$
\split
H^1(\Gamma; M) &= H^0(\overline{\Gamma}; H^1(\IAn; M))\oplus 
H^1(\overline{\Gamma}; M)\\
&= \Hom(\tensors, M)^{\overline{\Gamma}}\oplus 
H^1(\overline{\Gamma}; M).
\endsplit
$$
This completes the proof.\qed
\enddemo

When $M$ is a non-trivial irreducible $\bQ [\GL(\HZ)]$-module, as was proved 
by Borel [B], 
the cohomology group $H^*(\GL(\HZ); M)$ vanishes in a stable range, 
so that we obtain {\it stably}
$$
H^1(\An; M) = \Hom(\tensors, M)^{\GL(\HZ)}. \tag 6.7
$$
In the simplest case $M = H_\bQ$, we have 
$H^1(\An; H_\bQ) = \bQ$
for any {\it sufficient large} $n$. 
It is generated by the class $\h_1 = {\a_1}_*[\tau^\theta_1]$. 
On the other hand, Satoh [Sa1] used a direct method involved with a
presentation of the group $\An$ given by  Gersten [G] to prove that
$$
H^1(\An; \HZ) = \bZ
$$
for $n \geq 4$, and that it is generated by $\h_1$. 
Moreover he proved
$$
H^1(\Aut(F_3); \HZ) = \bZ/2\oplus\bZ.
$$

\beginsection 7. Decomposition of Cohomology Groups.\par

In this section we prove the cohomology class 
$\h_1 = {\a_1}_*[\tau^\theta_1] \in H^1(\An; H)$ 
gives a canonical decomposition of cohomology
groups of $\An$, and discuss the rational cohomology 
of $\An$ with trivial coefficients.\par
Recall the homomorphism $\iota: \Fn \to \An$,
$\gamma \mapsto (\iota(\gamma): \delta \mapsto
\gamma\delta\gamma\inv)$ in (2.9).
We have a group extension
$$
\Fn \overset\iota\to\to \An \overset\pi\to\to \On.
\tag 7.1
$$
\proclaim{Theorem 7.1} Suppose $1-n$ is invertible in the
ring $R$. Then we have a natural decomposition of the 
cohomology group
$$
H^*(\An; M) = H^*(\On; M)\oplus H^{*-1}(\On; H^*\otimes M)
$$
for any $R[\On]$-module $M$.
Especially,  $\pi^*: H^*(\On; M)\to H^*(\An;\allowmathbreak M)$ 
is an injection.
\endproclaim
\demo{Proof} 
We denote by $\{\ell_i\}^n_{i=1} \subset H^*$ the dual
basis of $\{X_i\}^n_{i=1} \subset H$. From (2.11) we have
$$
\tau^\theta_1(\iota(x_i)) = {\sum}^n_{j=1}\ell_j\otimes
(X_iX_j - X_jX_i)\in H^*\otimes H^{\otimes 2},
$$
and so $\a_1(\tau^\theta_1(\iota(x_i))) = (1-n)X_i$. Hence
$$
\iota^*\h_1 = (1-n)1_H \in H^1(\Fn; H) =\Hom(H, H). \tag 7.2
$$
The Lyndon-Hochschild-Serre spectral sequence of the
extension (7.1) induces an exact sequence
$$
\cdots \to H^p(\On; M)\overset{\pi^*}\to\to H^p(\An; M)
\overset{\pi_\sharp}\to\to H^{p-1}(\On; H^*\otimes M) \to\cdots.
\tag 7.3
$$
Here $\pi_\sharp$ is the Gysin map, and 
we have a natural isomorphism $H^1(\Fn; M) = H^*\otimes M$.
Consider the contraction map 
$$
C: H\otimes H^*\otimes M \to M, \quad
Y\otimes f \otimes m \mapsto f(Y)m.
$$ 
By (7.2) we have $\pi_\sharp(\h_1) = (1-n)1_H \in H^0(\On; H^*\otimes H)$, 
and so
$$
\pi_\sharp(C(\h_1\cup \pi^*u)) 
= (1_{H^*}\otimes C)(\pi_\sharp(\h_1)\cup u) 
= (1_{H^*}\otimes C)((1-n)1_H\cup u) 
= (1-n)u
$$
for any $u \in H^{p-1}(\On; H^*\otimes M)$. 
This implies the map
$$
H^{p-1}(\On; H^*\otimes M)\to H^{p}(\An; M), \quad
u \mapsto \frac{1}{1-n} C(k\cup \pi^*(u))
$$
is a right inverse of the map $\pi_\sharp$. 
Hence the exact sequence (7.3) splits. 
This completes the proof.\qed
\enddemo

We conclude the paper by discussing the rational cohomology 
of the group $\An$ with trivial coefficients. 
Similar results hold for the group $\On$ because 
$\pi^*: H^*(\On; \bQ) \to H^*(\An; \bQ)$ is injective 
by Theorem 7.1.\par
As was shown by Morita [Mo4], 
we obtain any of the Morita-Mumford classes on the mapping class group 
for a surface by contracting the coefficient of a power of the 
cohomology class $h_1$ in a suitable way involved 
with the intersection form $\HZ^{\otimes 2}\to \bZ$. 
Thus the class $h_1$ yields rich nontrivial classes  
in the rational cohomology of the mapping class group. 
For details, see also [KM][KM2].\par
For the group $\An$, in contrast, we have
\proclaim{Theorem 7.2}
$$
f_*({h_1}^{\otimes m}) = 0 \,\in \,H^m(\An; \bQ)
$$
for any $m \geq 1$ and any $\GL(\HZ)$-invariant linear form 
$f:(\tensors)^{\otimes m}\to \bQ$.
\endproclaim
\demo{Proof}
We denote by $M$ the $\GL(H_\bC)$-module 
$\Hom_\bC({H_\bC}^*\otimes\La^2H_\bC, \bC)$. 
Since the special linear group $\SL(\HZ)$ is 
Zariski dense in $\SL(H_\bC)$, we have 
$$
\Hom_\bZ((\tensors)^{\otimes m}, \bQ)^{\GL(\HZ)} \subset
(M^{\otimes m})^{\GL(\HZ)} \subset (M^{\otimes m})^{\SL(H_\bC)}.
\tag 7.4
$$
First we prove 
$$
(M^{\otimes m})^{\GL(\HZ)} = 0, \quad
\text{if \, $m \leq 2n-1$.}\tag 7.5
$$
Consider the tori
$$
\split
& \bT:=\{\diag(\zeta_1, \zeta_2, \dots, \zeta_n); \, 
\zeta_i \in \bC,\,\, \vert\zeta_i\vert = 1, \,\, 1 \leq i \leq n\}, 
\quad\text{and}\\
& \bT_0 :=\{\diag(\zeta_1, \zeta_2, \dots, \zeta_n)\in \bT; \, 
\zeta_1\zeta_2\cdots\zeta_n = 1\}\\
\endsplit
$$
in $\GL(H_\bC) = \GL_n(\bC)$. 
The representation ring $\R\bT$ is isomorphic to the Laurant polynomials
$$
\R\bT = \bZ[{t_1}^{\pm1}, {t_2}^{\pm1},\dots, {t_n}^{\pm1}].
$$
Here $t_i$ is the $1$-dimensional $\bC[\bT]$-module 
on which $\diag(\zeta_1, \zeta_2, \dots, \zeta_n)$ acts by 
multiplication by $\zeta_i$. Then $M$ is equal to 
$$
g(t) := (n-1){\sum}^n_{j=1}{t_j}\inv 
+ {\sum}_{i\neq j\neq k\neq i} t_i{t_j}\inv{t_k}\inv \in \R\bT
$$
as a $\bT$-module, which is homogeneous of degree $-1$. 
A monomial, or a $1$-dimensional representation, $h(t) \in \R\bT$ with
$\nu_i:= \deg_{t_i}h(t)\neq \nu_j:= \deg_{t_j}h(t)$ for some $i < j$ is
\underbar{not} 
$\bT_0$-invariant. In fact, $\diag(1, \dots, 1, \zeta, 1, 
\dots, 1, \zeta\inv, 1, \dots, 1)\in \bT_0$, where $\zeta$ is in the $i$-th
component, and $\zeta\inv$ in the $j$-th, 
acts on $h(t)$  by multiplication by $\zeta^{\nu_i- \nu_j}$.
Hence, if $n$ does not divide $m$, then $g(t)^m$ has no 
$\bT_0$-invariant part. 
Moreover, if $m=n$, $(M^{\otimes n})^{\bT_0}$ is equal to 
$c{t_1}\inv{t_2}\inv\cdots{t_n}\inv \in \R\bT$ for some $c \geq 0$. 
But ${t_1}\inv{t_2}\inv\cdots{t_n}\inv$ is \underbar{not} invariant 
under the action of $\diag(-1, 1, \dots, 1) \in \GL(\HZ)\cap \bT$. 
This proves (7.5).\par

Now, if $m \leq 2n-1$, the theorem follows from (7.5). 
On the other hand, as was established by Culler and 
Vogtmann [CV], the virtual cohomology dimension of $\An$
is $2n-2$. Hence, if $m \geq 2n-1$, then $H^m(\An; \bQ) = 0$. 
This completes the proof of the theorem.\qed
\enddemo

Recently Galatius [Ga] proved 
the rational reduced cohomology $\widetilde{H}^k(\An; \bQ)$ 
vanishes in a stable range, $n > 2k+1$. 
Hence our Theorem 7.2 for the stable range is an immediate 
consequence of Galatius' quite excellent result. 
On the other hand, instead of using his result, 
combining Theorem 7.2 with a result of Igusa [I], 
Theorem 8.5.3, p.333, and 
a computation of stable twisted cohomology of
$\GL(\HZ)$ by Borel [B], we can deduce that the homomorphism 
$\widetilde{H}^*(\An/[\IAn, \IAn]; \bQ) \to \widetilde{H}^*(\An; \bQ)$
induced by the natural projection vanishes in some stable range
$* \ll n$. Anyway we should remark that Theorem 7.2 holds also 
for the unstable range. In other words, the first Johnson map yields
no nontrivial rational cohomology classes on the automorphism group 
of the free group, $\An$, even in the unstable range, while it 
yields all the Morita-Mumford classes on the mapping class group 
for a surface. \par

%\newpage

\par
\vskip 15mm
\noindent
Department of Mathematical Sciences,
\newline
University of Tokyo
\newline
Tokyo, 153-8914 Japan 
\newline
e-mail address: kawazumi\@ms.u-tokyo.ac.jp
\par

\end